\documentclass[12pt]{article}
\usepackage{amssymb,latexsym,amsmath}
\usepackage{amsthm}
\usepackage{epsfig}
\usepackage{comment}
\usepackage{accents}
\usepackage{color}

\def\capa{{\rm cap \, }}
\def\dist{{\rm dist \, }}
\def\supp{{\rm supp \, }}

\def\Im{{\rm Im \, }}
\def\BB4#1{\textcolor{red}{#1}}

\newcommand{\C}{\mathbb{C}}

\newcommand{\N}{\mathbb{N}}
\newcommand{\R}{\mathbb{R}}

\newcommand{\MM}{\mathcal{M}}

\newcommand{\TT}{\mathcal{T}}

\renewcommand{\O}{\mathcal{O}}

\renewcommand{\Im}{{\rm Im} \,}

\newcommand{\Ker}{{\rm Ker} \,}

\newcommand{\conj}[1]{\left\{ {#1}\right\}}

\numberwithin{equation}{section}

\newtheorem{theorem}{Theorem}[section]
\newtheorem{lemma}[theorem]{Lemma}
\newtheorem{corollary}[theorem]{Corollary}
\newtheorem{proposition}[theorem]{Proposition}

\newtheorem{Definition}[theorem]{Definition}
\newenvironment{definition}{\begin{Definition}\rm}{\end{Definition}}
\newtheorem{Remark}[theorem]{Remark}
\newenvironment{remark}{\begin{Remark}\rm}{\end{Remark}}
\newtheorem{Example}[theorem]{Example}
\newenvironment{example}{\begin{Example}\rm}{\end{Example}}
\newtheorem{Problem}[theorem]{Problem}

\begin{document}
\title{Equilibrium problems for vector potentials with
semidefinite interaction
matrices and constrained masses}
\author{Bernhard Beckermann,
Valery Kalyagin, \\ Ana C. Matos and Franck Wielonsky}
\date{}
\maketitle

\begin{abstract}
We prove existence and uniqueness of a solution to the problem
of minimizing the logarithmic energy of vector potentials
associated to a $d$-tuple of positive measures supported on
closed subsets of the complex plane. The assumptions we make
on the interaction matrix are weaker than the usual ones and we
also let the masses of the measures vary in a compact subset
of $\R_+^d$. The solution is characterized in terms of
variational inequalities. Finally, we review a few examples taken
from the
recent literature that are related to our results.
\end{abstract}


\vspace{0.5cm}

\noindent {\bf AMS Classification (2010): } 31A15, 31A05, 30E10, 42C05, 41A28

\section{Introduction}
Vector equilibrium problems in logarithmic potential theory
have been studied for a few decades and have shown to be
crucial in
the investigation of many problems in
approximation theory, like those involving multiple orthogonal
polynomials (e.g. Hermite-Pad\'e approximants, in particular
Angelesco and Nikishin systems). This approach has been very
fruitful in the analysis of
numerous questions in numerical or applied mathematics, to
name only a few,
eigenvalue distribution of Toeplitz matrices, models in random
matrix theory, determinantal processes, e.g. non-intersecting
random paths. Vector equilibrium problems were first
considered in \cite{GR1,GR2}. The book \cite{niso} contains a
nice
introduction to the subject. Equilibrium problems on general
locally compact spaces are studied in \cite{Oht, Zori1, Zori2}.

We first introduce some notations. Let $\mu$ be a (positive) Borel
measure with closed support in $\C$ and set
\begin{equation}
U^\mu (z)=\int\log \frac{1}{|z-x|} d\mu (x),
\label{pot}
\end{equation}
for its logarithmic potential. Assume that $\mu$ has not too
much mass at infinity (in a sense to be specified later), so
that the above integral converges for $|z-x|$ large. Then, the
logarithmic potential
is a superharmonic function from $\C$ to $(-\infty,\infty]$,
and
the energy of $\mu$ is defined as
\[
I (\mu)=\iint_{}^{}\log \frac{1}{|x-y|} d\mu (x)d\mu (y)
=\int U^\mu(x)d\mu(x)
>-\infty.
\]
For a subset $\Sigma$ of $\C$, let
\begin{equation}
\label{def-MM}
\MM(\Sigma)=\{\mu \text{ Borel measure, of finite mass,
supported in }\Sigma, \text{ and } I(\mu)<\infty\},
\end{equation}
and
$$
\MM_t(\Sigma)=\{\mu\in\MM(\Sigma),~\|\mu\|=t\},
$$
where $\|\mu\|$ denotes the total mass of the measure $\mu$.
For two measures $\mu,\nu \in
\mathcal M(\Sigma)$, we define the so-called mutual energy
\begin{equation}
           I(\mu,\nu)=\iint \log\frac{1}{|x-y|} d\mu(x)d\nu(y).
           \label{mut-ener}
\end{equation}
Again, if $\mu$ and $\nu$ do not have too much mass at
infinity, this integral converges for $|x-y|$ large, and is
well-defined in $(-\infty,+\infty]$.

Throughout, we let
\begin{equation}
\label{Def_Delta}
\Delta= (\Delta_{1},\ldots,
\Delta_{d}),\qquad \cup_{i=1}^d \Delta_i\varsubsetneq\C,
\end{equation}
be a $d$-tuple of closed non polar sets of $\C$, i.e. of
positive logarithmic capacities
$$
\capa (\Delta_i)>0,\quad i=1,\ldots,d,
$$
and we define the cartesian products
\[
\MM^{d} (\Delta)=\MM(\Delta_{1})\times \cdots \times \MM (
\Delta_{d}),\quad
\MM^{d}_1 (\Delta)=\MM_1(\Delta_{1})\times \cdots \times
\MM_1(\Delta_{d}).
\]
Assume for the moment that the $\Delta_i$, $i=1,\ldots,d$, are
compact sets.
For two $d$-tuples of measures
\begin{equation*}
\mu= (\mu_{1},\ldots
,\mu_{d})^{t}\in \MM^{d} (\Delta),\quad \nu= (\nu_{1},\ldots
,\nu_{d})^{t}\in \MM^{d} (\Delta),\quad
\end{equation*}
we define the mutual energy of $\mu$ and $\nu$ as
\[
J (\mu,\nu)=\sum_{j=1}^{d}I (\mu_{j},\nu_{j}),
\]
which is finite. Actually, the compactness of the $\Delta_i$
entails that the mutual energy of two measures of finite
energies is also finite.

Let $C= (c_{i,j})$ be a real
symmetric positive definite matrix of order $d$, such that
\begin{equation}\label{compatNS}
\forall (i,j),\quad\text{if}\quad \Delta_{i}\cap
\Delta_{j}\neq\emptyset \quad\text{then}\quad c_{i,j}\geq 0.
\end{equation}
The energy of $\mu$ with respect to the interaction matrix $C$
is defined as
\[
J (\mu)= J(C\mu,\mu)=\sum_{i,j=1}^d c_{i,j} I(\mu_i,\mu_j).
\]
Note that, because of (\ref{compatNS}), $J(\mu)$ is always
well defined (even if some of the components of $\mu$ have
infinite
energies).
Now, the extremal problem is the following :
\\[\baselineskip]
\emph{find
\[
J^{*}=\inf \{J (\mu),\quad \mu\in \MM_1^d(\Delta)\},
\]
and characterize the extremal tuple of measures $\mu^{*}$ in
$\MM_1^{d} (\Delta)$, for which the infimum is
attained.}\\[\baselineskip]
As the sets $\Delta_i$ are assumed to be of positive
capacity,
a solution $\mu^*$ to this problem, with $J^{*}=J(\mu^*)<\infty$,
exists, and it is unique. The proof of existence is based on
the fact that the mutual energy (\ref{mut-ener}) is
lower semi-continuous, which implies together with
(\ref{compatNS}), that $\mu
\mapsto J(\mu)$ is also lower semi-continuous. Moreover, the
map is strictly convex on the
set $\MM_1^{d}(\Delta)$, from which uniqueness follows, see
\cite[Propositions 5.4.1 and 5.4.2]{niso}.

A characterization of the solution can be given via the
so-called equilibrium conditions. For
that, we introduce the partial potentials
\begin{equation*}
 U_{i}^{\mu}
(x)=\sum_{j=1}^{d}c_{i,j}U^{\mu_{j}} (x),\quad i=1,\ldots ,d,
\end{equation*}
where the scalar potentials $U^{\mu_{j}} (x)$ have been
defined in \eqref{pot}.
Then,  $d$-tuple of  measures $\mu$ solves the minimization problem if and
only if there exist constants $w_{i}$, such that, for
$i=1,\ldots ,d$,
\begin{align}\label{}
\label{char-sol1}
U_{i}^{\mu} (x)\geq w_{i},\quad & \text{quasi-everywhere on }
\Delta_{i},\\
\label{char-sol2}
U_{i}^{\mu} (x)\leq w_{i},\quad & \text{everywhere on }\supp
(\mu_{i}),
\end{align}
where quasi-everywhere means everywhere up to a set of
capacity zero. Proofs of these results can be found in
\cite[Chapter 5]{niso}.
\begin{remark}
\label{well_partial}
For some $x\in \C$ it may happen
that $U^{\mu_j}(x)=+\infty$ for several indices $j$. However,
the partial potential $U_i^{\mu}$
is well defined quasi-everywhere since positive measures of
finite mass and compact support have a finite potential
quasi-everywhere, see \cite[Theorem III.16]{Tsuji}.
\end{remark}
Regarding applications, it is also very useful to consider an additional external field in equilibrium
problems. The main reference for the study of equilibrium problems in presence of an external
fields is the book \cite{SaTo}.

Let $Q=(Q_j)_{j=1,...,d}$ be a vector of lower semi-continuous
functions,
$$
Q_j:\Delta_j \to (-\infty,\infty],\qquad j=1,\ldots,d,
$$
and define the weighted energy of a tuple of measures
$\mu\in\MM^d(\Delta)$ in the presence of the external field
$Q$ as
   \begin{equation}
    \label{def-J_Q}
          J_Q(\mu) = J(\mu) + 2 \sum_{j=1}^d \int Q_j
          d\mu_j .
   \end{equation}
For $\mu\in\MM^d(\Delta)$, we have mentioned that
$J(\mu)=J(C\mu,\mu)$ is finite. By lower-semicontinuity, each
$Q_j$ is bounded from below on $\Delta_j$, $j=1,\ldots,d$.
Hence, the integrals in \eqref{def-J_Q} are well-defined and
$J_Q(\mu)>-\infty$. It can also be checked that, in
$\MM_1^d(\Delta)$, there exists at least one measure $\mu$
with $J_Q(\mu)<\infty$, see the proof of Theorem \ref{existence} (i).

Then, the extremal problem of minimizing the weighted energies
\begin{equation}\label{min-JQ}
\{J_Q(\mu),\quad \mu\in\MM_1^d(\Delta)\},
\end{equation}
is solved by a unique d-tuple of measures
$\mu^{*}\in\MM_1^d(\Delta)$, with $J_Q(\mu^*)<\infty$, and it
is characterized by the existence of constants $w_{i}^{Q}$, such
that, for $i=1,\ldots ,d$,
\begin{align}\label{}
\label{char-solQ1} U_{i}^{\mu^*} (x)+Q_{i} (x)\geq
w_{i}^Q,\quad&
\text{quasi-everywhere on }\Delta_{i},\\
\label{char-solQ2} U_{i}^{\mu^*} (x)+Q_{i}
(x)\leq w_{i}^Q,\quad & \text{everywhere on } \supp (\mu_{i}).
\end{align}
For a proof in the scalar case $d=1$, we refer to \cite{GR2}
and \cite[Theorem I.1.3]{SaTo}.
The vector problem with
external fields is considered in \cite{GR2}, see also
\cite{FPLL}.

In the past few years, generalizations of the above vector
equilibrium problems have appeared repeatedly in the
literature. By generalizations, we mean that the assumptions on
the
interaction matrix or on the masses were relaxed in various
ways.
For instance, in \cite{AKWa,BB},
    one allows for sets which no longer satisfy the
compatibility condition \eqref{compatNS}, since some
$\Delta_j$ are intervals with a common endpoint. In
\cite{A,AL,AKWa},
one considers interaction matrices which are only positive
semidefinite. In these papers, the authors also minimize $J$ not over the set $\MM_1^d(\Delta)$ of
tuples of
probability measures but over the set $$
\MM^d_K(\Delta)=\{\mu=(\mu_1,...,\mu_d)^t \in \MM^d(\Delta),~
\|\mu\|=(\|\mu_{1}\|,\ldots,\|\mu_{d}\|)^t\in K\}, $$ where
$K$ is
a non-empty compact subset of the set $\R^d_{+}$ of $d$-tuples
of
non negative real numbers.
In addition, one considers in
\cite{AKWa,BB,CKL,DK,DelDui,duits_kuijlaars,SaTo,So} extremal
problems with not necessarily compact sets $\Delta_j$. In the
papers
\cite{AKWa,DK,DelDui,duits_kuijlaars,So},
a  solution satisfying the extremal
properties \eqref{char-sol1}--\eqref{char-sol2} or
\eqref{char-solQ1}--\eqref{char-solQ2} could be exhibited
directly through some algebraic equation hence settling
the
problem of existence of a minimizer.

The goal of this paper is to provide a more systematic
approach, by showing existence, uniqueness, and
characterization of the
extremal solution for a large class of generalized equilibrium
problems. At this point, the following simple examples are
instructive, since they show that some care has to be taken
when
weakening the assumptions of the minimization problem.

\begin{example}\label{example_condenser}
   Consider the data
   $$
     C = \left[\begin{array}{cc} 1 & -1 \\ -1 & 1
     \end{array}\right], \quad
      \Delta_1=[-1/2,0],\quad \Delta_2=[0,1/2],
   $$
where $C$ is positive semidefinite, and the problem of finding
the minimum $J^*$ of the corresponding energy
    $$J(\mu)=I(\mu_1-\mu_2)\geq 0,\qquad
   \mu=(\mu_1,\mu_2)^t \in \MM_1^2(\Delta).
    $$
It is known that the same problem on the pair of subsets
$\Delta_{1,n}=[-1/2,-1/n]$ and $\Delta_{2,n}=[1/n,1/2]$,
$n\geq1$, admits the minimal energy $J_n^*$ with
$$J_n^*=\frac{1}{\capa(\Delta_{1,n},\Delta_{2,n})}=\frac{2\pi
K(2/n)}{K'(2/n)},$$
where $\capa(\Delta_{1,n},\Delta_{2,n})$ denotes the capacity
of the condenser with plates $\Delta_{1,n}$ and $\Delta_{2,n}$.
The explicit value given in the second equality, in terms of the
complete and complementary elliptic integrals of the first
kind $K$ and $K'$, can be found in \cite{LL}, and
may also be derived from Example II.5.14 in
\cite[pp.133-134]{SaTo}. Since
$$K(k)=\frac{\pi}{2}+\O(k^2),\quad K'(k)=-\log k+\O(1),\qquad
\text{as }k\to 0,
    $$
we obtain by letting $n$ tend to infinity, that $J^*=0$.
However, this value cannot be reached by a couple of measures
$(\mu_1,\mu_2)$ of finite energy since $I(\mu_1-\mu_2)=0$
would imply $\mu_1=\mu_2$, see Lemma~\ref{J_pos} below.
\qed\end{example}

More generally, for a rank $1$ interaction matrix $C=y y^t$
with $y\in \{ -1,1 \}^d$, our vector equilibrium problem
corresponds
to the
electrostatics of a condenser with external field, see, e.g.,
\cite[Chapter~VIII]{SaTo}. Here one usually assumes disjoint
$\Delta_j$ in order to ensure existence and uniqueness of an
extremal tuple of measures, though, as we will see below, we may
somewhat relax this
condition.

Next, we present three simple examples where
existence of an extremal tuple of measures holds but not uniqueness.

\begin{example}\label{example0_not_uniqueness}
    Consider the data
    \[
        C=\left[\begin{array}{cc} 1 & 1 \\ 1 & 1
     \end{array}\right], \quad  \Delta_{1}=\Delta_{2}=[-1,1],
     \quad K=\{(x,y)\in\R^2,x+y=1,x\geq 0,y\geq 0 \},
\]
then $J (\mu_{1},\mu_{2})=I (\mu_{1}+\mu_{2})$ is minimal over
$\MM^2_K(\Delta)$ for any couples
$(x\omega_{[-1,1]},y\omega_{[-1,1]})$, $x+y=1$, where
$\omega_{[-1,1]}$ denotes the equilibrium
measure of $[-1,1]$.
 \qed
\end{example}

Here, one may show that $J$ is convex but not strictly convex
over
$\MM^2_K(\Delta)$.
Notice also that there is even not a unique minimizer over
$\MM^2_1(\Delta)$.

\begin{example}\label{example1_not_uniqueness}
    Consider the data
    \[
C=I_{2},\quad  \Delta_{1}=\Delta_{2}=[-1,1],
\quad K=\{(x,y)\in\R^2,x^{2}+y^{2}=1,x\geq 0,y\geq 0 \}.
\]
Then,
\[
J (\mu_{1},\mu_{2})=I (\mu_{1})+I (\mu_{2}),
\]
which is minimal when both measures $\mu_{1}$ and $\mu_{2}$
are multiples of the equilibrium measure $\omega_{[-1,1]}$ of
$[-1,1]$.
Hence, any couple $(x\omega_{[-1,1]},y\omega_{[-1,1]})$ with
$x^{2}+y^{2}=1$ belongs to $\MM_{K}^{2} (\Delta)$ and gives
the minimum value $\log 2$ of the energy $J$. \qed
\end{example}

For this example, it is not
difficult to show that $J$ is strictly convex over
$\MM^2(\Delta)$, but the non--uniqueness of the extremal tuple of measures comes from the lack of convexity of $K$. The next
example shows
that even convexity of $K$ does not allow to conclude.

\begin{example}\label{example2_not_uniqueness}
    Consider the data
    \[
        C=I_{2},\quad  \Delta_{1}=\Delta_{2}=[-4,4],
\quad K=\{(x,y)\in\R^2,x+y=1,x\geq 0,y\geq 0 \},
\]
then $J (\mu_{1},\mu_{2})=I (\mu_{1})+I (\mu_{2})$ is minimal
when both measures $\mu_{1}$ and $\mu_{2}$ are multiples of the
equilibrium measure $\omega_{[-4,4]}$ of $[-4,4]$, and in this case
$J(x\omega_{[-4,4]},y\omega_{[-4,4]})=(x^2+y^2) I(\omega_{[-4,4]})$.
Since $I(\omega_{[-4,4]})=-\log(2)<0$, we get
the
minimal value $-\log 2$ both for $(\omega_{[-4,4]},0)$ and
$(0,\omega_{[-4,4]})$ (and $J$ is no longer convex).
 \qed
\end{example}

In this work, we want to extend the afore-mentioned results
about the minimization of (\ref{min-JQ}) to the following
situation:\\[0.5\baselineskip]
(i) The sets $\Delta_i$, $i=1,\ldots,d$, are closed sets of
$\C$ (instead of compact sets).\\[0.5\baselineskip]
(ii) The interaction matrix $C\in \mathbb R^{d\times d}$, of
rank $r$ say, is positive semi-definite (instead of
definite).\\[0.5\baselineskip]
(iii) The compatibility condition \eqref{compatNS} is not
necessarily satisfied.
\\[0.5\baselineskip]
(iv) The minimization of $J_Q$ is performed over
$\MM_K^d(\Delta)$ instead of $\MM_1^d(\Delta)$.
\\[0.5\baselineskip]
To cope with the non-compactness of the sets $\Delta_i$ we need to add
to the defining properties of the set
$\MM(\Sigma)$, see \eqref{def-MM}, a growth condition at infinity.

Hence, from now on, the set $\MM(\Sigma)$ will consist of Borel measures $\mu$ of finite mass, supported on $\Sigma$, of finite energy, and such that
\begin{equation}\label{cond-mu}
\int\log(1+|t|)d\mu(t)<\infty.
\end{equation}
The set of $d$-tuples of measures $\MM_K^d(\Delta)$ is
redefined accordingly, i.e. we assume that condition
(\ref{cond-mu}) is
satisfied component-wise.

For a positive measure $\mu$ of finite mass, satisfying
(\ref{cond-mu}), we have
\begin{equation}
\label{pot-lower-bd}
U^\mu(z)\geq -\|\mu\|\log(1+|z|)-\int\log(1+|t|)d\mu(t)
>-\infty,\qquad z\in\C.
\end{equation}
The question raised in Remark \ref{well_partial} about the
well-definedness of partial potentials can be answered in the
same manner since the assertion given there still holds true
for measures in $\MM(\Sigma)$, see Lemma \ref{lem-tsuji}. For
two measures $\mu$ and $\nu$ of finite masses, satisfying
(\ref{cond-mu}), we have $$I(\mu,\nu)\geq
-\|\mu\|\int\log(1+|t|)d\nu(t)-\|\nu\|\int\log(1+|t|)d\mu(t)
>-\infty,
$$
and in particular $I(\mu)>-\infty$. Moreover, denoting by
$\widetilde\mu$ the normalized measure $\mu/\|\mu\|$ for a
non-zero $\mu\in\MM(\Sigma)$, it is known that the inequality
$$
I(\widetilde\mu-\widetilde\nu)\geq
0,\qquad\mu,\nu\in\MM(\Sigma),
$$
holds true, see Lemma \ref{J_pos}. In particular, we have $2I(\widetilde{\mu},\widetilde{\nu})\leq I
(\widetilde{\mu})+I(\widetilde{\nu})$, and since, by definition of
$\MM(\Sigma)$, the energies of $\mu$ and $\nu$ are finite, it
then follows that the mutual energy $I(\mu,\nu)$ is finite as
well. As a consequence, for $\mu\in\MM_K^d(\Delta)$, the energy $J(\mu)$ is always well defined
in $\R$.

For the external fields, we also need some growth condition at
infinity. Throughout, we assume that $Q=(Q_j)_j$ is a vector
of {\it{admissible}} functions,
in the sense\footnote{Compare with the slightly weaker growth condition at infinity given in \cite[Definition~I.1.1]{SaTo} for scalar extremal problems.} of \cite[Chapter VIII.1]{SaTo} :
\begin{definition}
\label{def-adm}
Let $\Sigma$ be a closed subset of $\C$ of positive capacity. A function $f: \Sigma
\to (-\infty,\infty]$ is said to be \emph{admissible} if it satisfies
the following three conditions:\\[0.5\baselineskip]
(i) $f$ is lower semi-continuous,\\[0.5\baselineskip]
(ii) $f$ is finite on a set of positive
capacity,\\[0.5\baselineskip]
(iii) $f(x)/\log|x|\to\infty$ as $|x|\to\infty$ (in case
$\Sigma$ is unbounded).
\end{definition}
In view of the preceding examples, we also have to add
assumptions\footnote{
In particular, Example~\ref{example0_not_uniqueness} tells us
that
the classical condition \eqref{compatNS} only ensures strict
convexity in case of invertible interaction matrices.
   }
linking the matrix of interaction $C$ to the topology of the
sets
$\Delta_j$. For the proof of the existence of an extremal
tuple of
measures we will assume that
\begin{equation}\label{assumpH2}
   \exists y\in\Im (C),\quad \forall (i,j),\quad
   \text{if}\quad
   \text{dist}(\Delta_i,\Delta_j)=0
   \quad\text{then}\quad y_{i}y_{j}> 0,
\end{equation}
whereas, for uniqueness, we will also impose that, for
any subset of indices $I\subset \{ 1,2,\ldots,d \}$, different
from a singleton,
\begin{equation}\label{assumpH1}
  \mbox{if the columns $(C_i)_{i\in I}$ of $C$ are linearly
  dependent then} \quad \capa\Bigl( \bigcap_{i \in I} \Delta_i
  \Bigr) = 0.
\end{equation}
Notice that both conditions (\ref{assumpH2}) and
(\ref{assumpH1})
are trivially true
for positive definite interaction matrices $C$ (for condition
(\ref{assumpH2}) take
$y=(1,...,1)^t$). Such interaction matrices appear e.g. when
studying the asymptotic behavior of Angelesco or Nikishin
systems in approximation theory.

It is instructive to have a closer look at vector equilibrium
problems corresponding to condensers, namely with interaction
matrices $C=yy^t$ of rank $1$, $y \in \{ -1,1 \}^d$. In this
case,
(\ref{assumpH2}) is equivalent to \eqref{compatNS}, it tells
us
that any two plates $\Delta_j$ with charges of opposite sign
have
positive distance, and (\ref{assumpH1}) requires in addition
that
any two plates $\Delta_j$ with charges of the same sign have
an intersection of capacity zero. Finally, notice that condition
(\ref{assumpH2}) fails to hold for
Example~\ref{example_condenser}, whereas condition
(\ref{assumpH1}) fails to hold for
Example~\ref{example0_not_uniqueness}. For the other two
examples,
conditions (\ref{assumpH2}) and (\ref{assumpH1}) hold,
indicating
that there should be additional restrictions on the set $K$.

We now state the two main results of our paper. The first result shows, under assumption (\ref{assumpH2}), the existence of a solution to our mimimization problem.
\begin{theorem}\label{existence}
Consider some nonempty compact set $K\subset \mathbb R_+^d$,
and assume that the positive semidefinite interaction matrix
$C$ satisfies \eqref{assumpH2}. Let
\begin{equation}
\label{min-pb}
J_{Q}^* : =\inf \{J_{Q} (\mu),~ \mu\in \MM^{d}_{K}
(\Delta)\}.
\end{equation}
Then, the following assertions hold. \\
(a) $J_Q^*$ is finite.\\
(b) There exists a $d$-tuple of measures
$\mu^{*}\in\MM^{d}_{K} (\Delta)$,
such that
$J_{Q} (\mu^{*}) =J_{Q}^*$.
\end{theorem}
Our second result is about uniqueness of a minimizer of the
extremal problem \eqref{min-pb}, and about its characterization by equilibrium
conditions, the so-called Euler--Lagrange inequalities. Here
we restrict ourselves to measures $\mu$ whose vector of masses
$(\|\mu_{1}\|,\ldots ,\|\mu_{d}\|)$ lies in a non--empty
compact polyhedron $K$ of $\R_+^d$.
\begin{theorem}\label{uniqueness} Assume that the positive
semidefinite interaction matrix $C$ satisfies the assumptions
\eqref{assumpH2} and \eqref{assumpH1}, and that the set of
masses
$K$ consists of a non--empty compact polyhedron of the form
\begin{equation}
\label{def-K}
   K = \{x\in\R^d_+, Ax= a \},
\end{equation}
with $A\in \mathbb R^{m \times d}$ and $a\in \mathbb R^m$,
where we suppose in addition that
\begin{equation} \label{image_AC}
          \Ker (A)\subset \Ker (C).
\end{equation}
Then, the following assertions hold true,\\
(a) There exists a unique $d$-tuple of measures $\mu^{*}\in\MM^{d}_{K}
(\Delta)$, of finite energy $J_{Q} (\mu^{*})<\infty$, such
that\[
J_{Q} (\mu^{*}) =\inf \{J_{Q} (\mu),\quad \mu\in \MM^{d}_{K}
(\Delta)\}.
\]
(b) The $d$-tuple of measures $$\mu=(\mu_1,...,\mu_d)\in
   \MM^d_K(\Delta)$$
   is the minimizer
of $J_Q$ over $\MM^d_K(\Delta)$ if and only if there exists $F
\in \mathbb R^m$ such that, for $i=1,\ldots,d$,
   \begin{align}
\label{char-solQC1} U_{i}^{\mu} (x)+Q_{i} (x)\geq (A^t
F)_i\quad&
         \text{quasi--everywhere on $\Delta_{i}$,}\\
         \label{char-solQC2} U_{i}^{\mu} (x)+Q_{i}
              (x)\leq (A^t F)_i\quad &
         \text{$\mu_i$--almost everywhere on $\Delta_{i}$.}
\end{align}
\end{theorem}
\begin{remark}
Notice that Theorem~\ref{uniqueness} includes the particular
case
$A=I_d$ of a singleton $K$, where we prescribe the mass of all
components of our tuple of measures.
Non--singleton $K$ of the form \eqref{def-K} have been
considered
first in \cite{A,AL,AKWa}, where the authors impose equality
in \eqref{image_AC}. From Example~\ref{example2_not_uniqueness}
we learn that in general the condition \eqref{image_AC} cannot be
dropped for establishing uniqueness. \qed
\end{remark}
As said before, in case of invertible $C$, all our (somehow technical) assumptions are trivially true for any configuration of sets $\Delta_j$ as in \eqref{Def_Delta}.
\begin{corollary}\label{coro}
  In case of a symmetric positive definite interaction matrix $C$ and a singleton $K=\{ a \}$, there exists one and only one minimizer
  of $J_Q$ over $\mathcal M_K^d(\Delta)$, which is characterized by the equilibrium conditions \eqref{char-solQC1} and \eqref{char-solQC2} for $A=I_d$.
\end{corollary}
\begin{example}\label{example_condenser2}
  Let
  $$
        C=\left[\begin{array}{cc} 2 & -1 \\ -1 & 2
     \end{array}\right], \quad \Delta_2 \subset \Delta_1 \subset \mathbb C,
      \quad K=\{ (a_1,a_2) \}\subset \mathbb R_+^2, \quad
      a_2\leq 2 \, a_1
  $$
  then, according to \eqref{char-solQC1} and \eqref{char-solQC2}, the couple of measures
  $$
        \mu_1 = (a_1-\frac{a_2}{2}) \omega_{\Delta_1} + \frac{a_2}{2} \omega_{\Delta_2} , \quad
        \mu_2 = a_2 \omega_{\Delta_2} ,
  $$
  minimizes $J$ over $\mathcal M_K^d(\Delta)$. As in Example~\ref{example_condenser}, we can give an electrostatic interpretation in terms of a condenser with two plates $\Delta_1$ and $\Delta_2$ of opposite charge. However, here the Nikishin interaction matrix $C$ translates some di-electric medium where particles of equal charge have stronger interaction than those of opposite sign. We observe the somehow surprising fact that there exists a unique  electrostatic equilibrium even if the two plates overlap. Notice that Nikishin systems with touching intervals $\Delta_j$ have been considered before in the literature without addressing this issue.
  \qed
\end{example}

The remainder of the paper is organized as follows. In Section
\ref{Prelim}, we gather several preliminary results that are
needed in the sequel. In Section \ref{Existence}, we give the proof of Theorem \ref{existence}. We also
derive,
under an additional condition, that the components of the
solution
have compact supports. In Section \ref{Unique}, we give the proof of Theorem \ref{uniqueness}.
In Section
\ref{Examples}, we review a few examples taken from the
literature that are
related to our results. Some open questions are discussed in
Section \ref{Conclusion}.

\section{Preliminary Results}
\label{Prelim}
Let us first recall the important fact that the energy of a
signed measure of mass 0 is non negative.
    \begin{lemma}
    \label{J_pos}
    Let $\mu,\nu\in\MM(\C)$ with $\|\mu\|=\|\nu\|$. Then
    $$I(\mu-\nu)\geq 0,$$
    and $I(\mu-\nu)=0$ if and only if $\mu=\nu$.
    \end{lemma}
\begin{proof}
See \cite[Lemma I.1.8]{SaTo} for measures $\mu,\nu$ with
compact support, and \cite[Theorem 2.5]{CKL} for the unbounded case, see also
\cite[Theorem 4.1]{Sim}.
\end{proof}

We proceed with a few results which are well-known when the
supports of the measures are compact, but for which we have
not always found references in the non compact case. We defer the
proofs
of these results to the Appendix.
    \begin{lemma}
    \label{mes_zero}
Let $\mu\in\MM(\C)$ . Then $\mu(E)=0$ for every Borel polar
set $E$.
    \end{lemma}
\begin{proof}
 See \cite[Remark~I.1.7]{SaTo} for $\supp(\mu)$ compact, and
   Appendix~\ref{Proofs} for the general case.
\end{proof}
\begin{lemma} \label{lem-tsuji}
   Let $\mu$ be a positive measure of finite mass, satisfying
(\ref{cond-mu}). Then, the potential $U^{\mu}(z)$ can be
$+\infty$
   only on a Borel  set of capacity 0.
\end{lemma}
\begin{proof}
It is well-known that the assertion holds true for any super-harmonic function on $\C$, not
identically $+\infty$, see \cite[Theorem 3.5.1]{ran}. In particular, it holds true for the potential $U^
{\mu}(z)$.
\end{proof}

Throughout we will use weak convergence of Borel measures. Let $(\mu_n)_n$ be a bounded sequence of Borel
measures on $\C$,
$$\|\mu_n\|\leq c<\infty,\qquad n\in\N.
$$
We recall that the sequence $\mu_n$ tends weakly to a measure
$\mu$, as $n\to\infty$, if
\begin{equation}\label{def-weak}
\int fd\mu_n\to\int fd\mu,
\end{equation}
for every bounded, continuous, real-valued function $f$ on
$\C$. In the literature, the notion of vague convergence is
also used,
where it is assumed that (\ref{def-weak}) holds true only for
continuous function $f$ on $\C$ with compact support. Clearly,
vague convergence is
weaker than the weak convergence. For example, the
sequence $\delta_n$ of Dirac measures at $x=n$ converges
vaguely to 0, although it does not converge weakly.
For some comments on these two different notions of
convergence of measures, one may have a look to
\cite[pp.134--137]{Deift}.
\begin{lemma}\label{weak-lower-sc}
Assume that the bounded sequence $\mu_n$ tends weakly to $\mu$, and let $Q$ be a lower bounded, lower
semi-continuous function on $\C$. Then
    $$\int Qd\mu\leq\liminf_{n\to\infty}\int Qd\mu_n.$$
    \end{lemma}
\begin{proof}
 See \cite[Theorem 0.1.4]{SaTo} for $\mu_n$ all supported in
   a compact set, and
   Appendix~\ref{Proofs} for the general case.
\end{proof}
    \begin{definition}
    \label{def-tight}
A bounded sequence of measures $(\mu_n)_{n\geq 0}$ in
$\MM(\C)$ is said to be \\
(i) tight if:
        \begin{equation*}
\forall \epsilon>0, \quad\exists \text{ compact set
}K\subset\C,\quad\forall n\in\N,\quad
\int_{\C\setminus K}d\mu_n(t)\leq\epsilon,
\end{equation*}
(ii) log-tight if:
\begin{equation}\label{log-tight}
    \forall \epsilon>0, \quad\exists
    \text{ compact set }K\subset\C, \quad\forall n\in\N,\quad
\int_{\C\setminus K}\log(1+|t|)d\mu_n(t)\leq\epsilon.
\end{equation}
\end{definition}
The notion of tightness of a bounded set of measures is
classical, see e.g. \cite{Bill}. The notion of log-tightness
is slightly stronger. Note that, from assumption \eqref{cond-mu},
each individual measure $\mu\in\MM(\C)$ satisfies inequality
(\ref{log-tight}).
Here, for log-tightness of a sequence, we ask this condition
to be satisfied uniformly with respect to $n$.

\begin{theorem}[Prohorov]\label{Prohorov}
Let $(\mu_n)_{n\geq 0}$ be a tight sequence of
probability measures on $\C$. Then, there is a subsequence of
$(\mu_n)_{n\geq 0}$ which is weakly convergent to a
probability measure on $\C$.
\end{theorem}
\begin{proof}
See Helly's selection theorem \cite[Theorem 0.1.3]{SaTo} for
$\mu_n$ all supported in
some compact set, \cite[Theorem 5.1]{Bill} in a general metric
space, and
   \cite[Theorem 9.3.3]{Dudley} for the special case of the
   euclidean space $\R^k$.
\end{proof}
\begin{remark}
The Prohorov theorem is actually stronger than Theorem \ref{Prohorov}, in that it also states, in the converse direction, that a weakly convergent sequence of measures is tight.
\end{remark}
    \begin{lemma}
    \label{scalar_lem_sc}
Let $(\mu_n)_{n\geq 0}$ and $(\nu_n)_{n\geq 0}$ be bounded
log-tight sequences of measures in $\MM(\C)$.
Assume $\mu$ and $\nu$ are two Borel measures such that $\mu_n\to\mu$ and $\nu_n\to\nu$ in
the weak
topology. Then
\begin{equation}
\label{I-lsc}
I(\mu,\nu)\leq\liminf_{n\to\infty}I(\mu_n,\nu_n).
\end{equation}
\end{lemma}
\begin{proof}
   See \cite[Theorem~5.2.1]{niso}
   for all $\mu_n,\nu_n$ supported in
   some compact set, and
   Appendix~\ref{Proofs} for the general case.
\end{proof}

Let us proceed with establishing four propositions, among
which we prove
the positiveness of $J$, the lower semi-continuity of $J_Q$,
and an inequality relating the weighted energy $J_Q(\mu)$ with
the
scalar energies of the components of $\mu$.

Throughout, we write the positive
semidefinite matrix $C$ of rank $r$ as
  a full rank factorization of the form
  \begin{equation}\label{def-C}
C=B^{t}B,\quad B\text{ matrix of dimensions }(r,d),~ r\leq
d,\text{ of rank } r.
  \end{equation}
Such a factorization is obtained, e.g., from the Jordan decomposition of $C$ by recalling that there exists an orthonormal basis of eigenvectors of $C$.  First, we generalize Lemma \ref{J_pos} to our vector setting.
\begin{proposition}\label{lem_bounded}
Let $\mu,\nu\in\MM_K^d(\Delta)$ with tuples of masses
verifying $B\|\mu\|=B\|\nu\|$. Then,
   \begin{equation} \label{positivity}
       J(\mu-\nu)\geq 0 .
   \end{equation}
Moreover, if condition \eqref{assumpH1} holds true, then
   \begin{equation}
        \label{positivity2}
J(\mu-\nu)= 0 \quad \mbox{if and only if}\quad \mu=\nu.
\end{equation}
\end{proposition}
\begin{proof}
  Let $\lambda = B(\mu-\nu)$. Then, we may write
   $$
        J(\mu-\nu) = J (\lambda,\lambda)=\sum_{j=1}^r
        I(\lambda_j).
   $$
By assumption, each component $\lambda_j$ of $\lambda$, with
Hahn decomposition
$\lambda_j=\lambda_{j,+}-\lambda_{j,-}$, is a signed measure
of mass $0$, whose absolute value $\lambda_{j,+}+\lambda_{j,-}$
is of
finite energy. Hence Lemma \ref{J_pos} applies, showing that
each
$I(\lambda_j)$ is non negative, so that \eqref{positivity}
holds true.

   We also know from Lemma \ref{J_pos} that $J(\mu-\nu)>0$ if
$\lambda_j\neq 0$ for at least one index $j$. Hence, to
establish
\eqref{positivity2} it only remains to show that $\mu\neq
\nu$ implies $\lambda\neq 0$. This property is trivial for
positive definite $C$ and thus invertible $B$. In our setting
with
semidefinite $C$, we will need assumption \eqref{assumpH1}.

Assuming $\mu -\nu \neq 0$, we deduce that there exists an
index $i_0$ and a Borel set $N$ such that
$(\mu_{i_0}-\nu_{i_0})(N)\neq 0$. Now, we consider the
partition$$
\bigcup_{j=1,\ldots,d}\Delta_j =
\bigcup_{I\subset\{1,\ldots,d\},~I\neq\emptyset}E_I,\quad
E_I=\left(\bigcap_{i\in
I}\Delta_i\right)\cap\left(\bigcap_{i\notin
I}\Delta_i^c\right),$$
where some of the $E_I$ may be empty sets. This induces a
partition of $N$,
   $$
N = \bigcup_{I \subset \{1,...,d \},~I\neq\emptyset} N_I,
\quad N_I=N\cap E_I,
   $$
so that $(\mu_{i_0}-\nu_{i_0})(N)=\sum_I
(\mu_{i_0}-\nu_{i_0})(N_I)$.
Therefore there exists a subset $I\subset \{ 1,...,d \}$ such
that
   $(\mu_{i_0}-\nu_{i_0})(N_I)\neq 0$, and
   \begin{equation}\label{mu-nu}
\forall i\notin I,\quad (\mu_{i}-\nu_{i}) (N_I) =
(\mu_{i}-\nu_{i}) (N_I\setminus \Delta_i) = 0,
   \end{equation}
     since $\supp(\mu_i-\nu_i)\subset \Delta_i$.
Note also that either $\mu_{i_0}(N_I)$ or $\nu_{i_0}(N_I)$ is
nonzero, so that $N_I$ is of positive capacity by Lemma
\ref{mes_zero}.
  Denote by $\widetilde{B}$, and $\widetilde C$,
the submatrix of $B$, and of $C$, respectively, obtained from
selecting
the columns of indices belonging to $I$. Since with $N_I$
also $\bigcap_{i\in I} \Delta_i$ has positive capacity, we
obtain
from
   condition \eqref{assumpH1} that $\widetilde C$ and
   thus $\widetilde B$ has full column rank.
By \eqref{mu-nu}, the relation $B (\mu-\nu) (N_I)=
\lambda(N_I)$ simplifies to $\widetilde{B}(\mu-\nu)_{i\in I}
(N_I)=\lambda(N_I)$, which cannot be zero. Thus $\lambda \neq
0$.
\end{proof}
As in the classical case, the main ingredient in the proof of
Theorem \ref{existence} will be the lower semi-continuity of
the functional $J_{Q}$.
We note that the proof does not
use the compatibility condition
\eqref{compatNS}.
\begin{proposition}\label{lem_semicontinuity}
Let $(\mu^{(n)})_{n\geq 0}$ be a sequence of $d$-tuples of measures in
$\MM_{K}^{d}(\Delta)$ which is log-tight (in the
component--wise sense)
and assume that $\mu^{(n)}$ tends to a $d$-tuple of measures $\mu\in\MM_K^d(\Delta)$, again
component--wise, as $n\to\infty$, in the weak topology.
Then
$$J_Q(\mu)\leq\liminf_{n\to\infty}J_Q(\mu^{(n)}).$$
\end{proposition}
\begin{proof}
We first show the asserted inequality for the map $\mu \mapsto
J(\mu)$. For that, we will use convolution of scalar finite
Borel
measures $\mu$ and $\nu$, which, for a Borel set $B\subset\C$, is defined as follows,
   $$
        (\mu * \nu)(B) = \int
        \nu(B-t) d\mu(t) =\int
        \mu(B-t) d\nu(t).
   $$
The convolution $\mu * \nu$ is a positive measure such that
$$
\supp(\mu*\nu)\subset\supp(\mu)+\supp (\nu),
\quad(\mu*\nu)(\C)=\mu
(\C)\nu (\C).
$$
From
$$(\mu * \nu)(B)=(\mu\times\nu)(\{(x,y),~x+y\in B\},
$$
it is easy to see that convolution is a commutative and
associative operation.
We will also use the convolution of a function $h$ with
a measure $\mu$,
\[
h*\mu (z)=\int_{}^{}h (z-t)d\mu (t),
\]
so that the potential $U^{\mu}$ coincides with the convolution
$-\log|\cdot|*\mu$.

Let $\lambda_N$ be the equilibrium measure of the
   circle centered at 0 of radius $e^{-N}$. Its potential is easily computed:
\[
U^{\lambda_N}(x)=\min\left(N,
   \log\frac{1}{|x|}\right),
\]
see e.g. \cite[Example 0.5.7]{SaTo}.
It is a continuous function tending pointwise to $\log
(1/|x|)$, $x\neq 0$, as
$N$ tends to $\infty$. Then, by associativity and
commutativity of the
convolution, we get
   \begin{equation*}
         U^{\mu * \lambda_N}(z)= -\log|\cdot|*
(\mu*\lambda_N)(z)= ( -\log|\cdot|*\lambda_N)*\mu (z)= \int
         U^{\lambda_N}(z-x)d\mu(x),
\end{equation*}
and for the mutual energies, we have
\begin{align}
         I(\mu * \lambda_N,\nu) & = \int U^{\lambda_N}(x-y)
         d\mu(x) d\nu(y),
         \\ \notag
         I(\mu * \lambda_N,\nu* \lambda_N)
& =\int U^{\lambda_N}(x-y) d\mu(x) d (\nu*\lambda_N)(y)
=\int
(U^{\lambda_N}*(\nu*\lambda_N))(x)d\mu(x) \\
& =\int U^{\lambda_N*\lambda_N}(x-y) d\mu(x) d\nu(y).
   \end{align}
From the definition of $U^{\lambda_N}$ follows that
$I(\mu * \lambda_N,\nu* \lambda_N)\leq I(\mu,\nu)$.
In particular, $I(\mu * \lambda_N)<\infty$ if $I(\mu)<\infty$.
Moreover,
\begin{align*}
\int\log(1+|x|)d(\mu * \lambda_N) & =
\iint\log(1+|x+y|)d\mu(x)d\lambda_N(y)\leq
\int\log(1+e^{-N}+|x|)d\mu(x)\\
& \leq  \log(1+e^{-N})\|\mu\|+\int\log(1+|x|)d\mu(x)<\infty.
\end{align*}
Hence, for any closed subset $\Sigma$ of $\C$, the measure $\mu*
\lambda_N$ lies in $\MM(\Sigma+D (0,e^{-N}))$ if
$\mu\in\MM(\Sigma)$.

Now, consider a log-tight sequence $\mu^{(n)}\in\MM_K^d(\Delta)$ such that
   $$
\mu^{(n)} \to \mu \in \mathcal
M_{K}^{d}(\Delta),
   $$
in the weak sense. Let $N$ be given. From the above
remarks,
   the $d$-tuple of measures $\mu^{(n)} * \lambda_N$, where
the convolution is taken componentwise, belongs to
$\mathcal M_{K}^{d}(\Delta+D (0,e^{-N}))$, and the masses of
$\mu^{(n)}$ and
$\mu^{(n)} * \lambda_N$ are the same. Thus, from
(\ref{positivity}), we get
   \begin{equation*}
         J(\mu^{(n)}-\mu^{(n)}* \lambda_{N})\geq 0,
                \end{equation*}
                or equivalently,
                \begin{equation*}
            J(\mu^{(n)})
                \geq         \sum_{i,j=1}^d c_{i,j} \Bigl(
         I(\mu_i^{(n)},\mu_j^{(n)}* \lambda_N)
         + I(\mu_i^{(n)}* \lambda_N,\mu_j^{(n)})
         - I(\mu_i^{(n)}* \lambda_N,\mu_j^{(n)}* \lambda_N)
         \Bigr).
   \end{equation*}
Let us consider the first energy in the right-hand side of the
above inequality. Since $\mu_i^{(n)}* \lambda_N$ is a
log-tight family which tends weakly to $\mu_i* \lambda_N$,
Lemma
\ref{scalar_lem_sc} tells us that,
\begin{equation}\label{conv-inf-I}
\liminf_{n\to\infty}I(\mu_i^{(n)},\mu_j^{(n)}* \lambda_N)\geq
I(\mu_i,\mu_j* \lambda_N).
\end{equation}
Actually, we have more. Indeed, redoing the proof of Lemma
\ref{scalar_lem_sc} with the kernel $U^{\lambda_N}(x-y)$
instead of
$\log(|x-y|^{-1})$, we now get an integrand in the first
integral of (\ref{eq1}) which is bounded and continuous.
Hence, (\ref{conv-inf-I}) can be strengthened to
\begin{equation*}
\lim_{n\to\infty}I(\mu_i^{(n)},\mu_j^{(n)}* \lambda_N)=
I(\mu_i,\mu_j* \lambda_N).
\end{equation*}
 The limits
$$\lim_{n\to\infty}I(\mu_i^{(n)}* \lambda_N,\mu_j^{(n)})=
I(\mu_i* \lambda_N,\mu_j),\quad
\lim_{n\to\infty} I(\mu_i^{(n)}* \lambda_N,\mu_j^{(n)}*
\lambda_N)=
I(\mu_i* \lambda_N,\mu_j* \lambda_N)
                $$
                are proven in the same way.
Consequently, we obtain that
\[
\liminf_{n \to \infty} J(\mu^{(n)})\geq
         \sum_{i,j=1}^d c_{i,j} \Bigl(
         I(\mu_i,\mu_j * \lambda_N)
         + I(\mu_i* \lambda_N,\mu_j)
         - I(\mu_i* \lambda_N,\mu_j* \lambda_N)
         \Bigr),
\]
where the right-hand side has a sense since we assume that the limit measure $\mu\in\MM_K^d
(\Delta)$ (all its components have finite energy).
Finally, both potentials $U^{\lambda_N}$ and
$U^{\lambda_N*\lambda_N}$ tend pointwise to
$\log (1/|x|)$ for $x\neq 0$, as $N\to \infty$. They
are dominated by $|\log (1/|x|)|$, and moreover,
\begin{multline*}
\iint\left|\log\frac{1}{|x-y|}\right|d\mu_{i}(x)d\mu_{j}(y)\\
\leq I(\mu_{i},\mu_{j})+
2\|\mu_{j}\|\int\log(1+|x|)d\mu_{i}(x)+2\|\mu_{i}\|\int\log(1+|x|)d\mu_{j}(x),
\end{multline*}
which is finite because
$I(\mu_i,\mu_j)$ is and we have (\ref{cond-mu}).
Hence, from
the dominated convergence theorem, we get
   \begin{equation*}
      \lim_{N\to \infty} I(\mu_i,\mu_j * \lambda_N)
      =I(\mu_i,\mu_j ), \quad
      \lim_{N\to \infty} I(\mu_i* \lambda_N,\mu_j * \lambda_N)
      =I(\mu_i,\mu_j ),
   \end{equation*}
   implying that
   $$
         \liminf_{n \to \infty} J(\mu^{(n)})\geq
                  J(\mu).
   $$
Since the external fields $Q_j$ are lower semi-continuous and
lower bounded, the fact that
\[
\liminf_{n \to \infty}\int Q_{j}d\mu_j^{(n)}\geq \int
Q_{j}d\mu_j,\quad j=1,\ldots ,d,
\]
follows from Lemma \ref{weak-lower-sc}.
\end{proof}
The aim of the next proposition is to show an inequality which
will be used in the proof of Proposition \ref{ineqJ-I}. It
asserts
that the scalar energy of a linear combination $\sum_j y_j\mu_j$ of
bounded measures $\mu_j$ in $\MM(\C)$, with given coefficients
$y_j$, is lower bounded, independently of the $\mu_j$, as soon
as
it is weighted by a multiple $\gamma Q$ of the external field,
with $\gamma$ an arbitrary small positive number.
Such a result is needed only to cope with unbounded $\Delta_j$
since, for compact $\Delta_j$, it is not difficult to derive a lower bound for the energy of a signed
measure which does not
involve external fields.
\begin{proposition}\label{lem_ineqQ-I}
Let $Q=(Q_1,\ldots,Q_d)^t$ be an admissible external field and
let $y=(y_1,\ldots,y_d)^t$ be a given vector in $\R^d$. Then,
\begin{equation}
\label{ineqQ-I}
\forall \gamma>0,\quad\exists \Gamma\in\R,\quad\forall \mu
\in \mathcal
M_{K}^{d}(\Delta),\quad \Gamma\leq I(y^t\mu)+\gamma
\int Q^td\mu.
\end{equation}
\end{proposition}
\begin{proof}
Since the union $\Sigma$ of the sets $\Delta_i$,
$i=1,\ldots,d$, is different from $\C$, recall (\ref{Def_Delta}), there exist some
$z_0\in\C$ and some
$r<1$ say, such that the disk
$D(z_0,2r)$ does not intersect $\Sigma$. Let $\omega_D$ be the
equilibrium measure of the disk $D=D(z_0,r)$ and let
$$\tau=\lambda-\lambda(\C)\omega_D,
$$
where $\lambda$ denotes the scalar signed measure $y^t\mu$. Since
$I(\lambda)$ is finite, $I(\tau)$ is finite as well and Lemma
\ref{J_pos} applies : $I(\tau)\geq 0$, or equivalently
\begin{align}\notag
I(\lambda) & \geq
2\lambda(\C)I(\lambda,\omega_D)+\lambda(\C)^2\log(r)\\
\label{ineqa}
& = 2\lambda(\C)\sum_{j=1}^dy_jI(\mu_j,\omega_D)+
\lambda(\C)^2\log(r).
\end{align}
All the mutual energies $I(\mu_j,\omega_D)$ can be bounded
above:
\begin{equation}\label{ineqb}
I(\mu_j,\omega_D)=\iint\log\frac{1}{|z-t|}d\mu_jd\omega_D\leq
\log\left(\frac1r\right)\|\mu_j\|\leq\log\left(\frac1r\right)M_j(K),
\end{equation}
with $M_j(K)=\sup_{\mu\in\MM_K^d(\Delta)}\|\mu_j\|$.
Moreover, the $I(\mu_j,\omega_D)$ can also be lower bounded.
First, note that, in view of the third condition of
admissibility in Definition \ref{def-adm} and the fact that
$Q_j$ is lower bounded on compact sets, we have
$$\forall\gamma_j>0,\quad\exists\Gamma_j\in\R,\quad\forall
z\in\Delta_j,\quad
\log(1+|z|)\leq\gamma_j Q_j(z)+\Gamma_j.$$
Then,
\begin{align}\notag
-I(\mu_j,\omega_D) & \leq\int\log(1+|z|)d\mu_j(z)+
\|\mu_j\|\int\log(1+|t|)d\omega_D(t)\\\notag
& \leq
\gamma_j\int Q_jd\mu_j+\Gamma_j\|\mu_j\|
+\|\mu_j\|\sup_{t\in D}\left(\log(1+|t|)\right)\\\label{ineqc}
& \leq \gamma_j\int Q_jd\mu_j+M_j(K)\left(\Gamma_j
+\sup_{t\in D}\left(\log(1+|t|)\right)\right),
\end{align}
and the proposition follows from plugging inequalities
(\ref{ineqb}) or (\ref{ineqc}) into (\ref{ineqa}), according
to the sign of
$\lambda(\C)y_j$, and noting that $\lambda(\C)$ is bounded
both above and below independently of $\mu$.
\end{proof}
Next, we show that the weighted energy of a tuple of measures
$\mu\in\MM_K^d(\Delta)$ dominates
the energies of its components. This result requires
the condition \eqref{assumpH2}.
\begin{proposition}\label{ineqJ-I}
Assume that the $d$-tuple of closed sets $\Delta$ and the
interaction matrix $C$ satisfy \eqref{assumpH2}. Then, there
exist positive constants $a_0$ and $a_{1}$ such that
\begin{equation}\label{ineqJI}
       \forall \mu= (\mu_{1},\ldots ,\mu_{d})^{t} \in \mathcal
M_{K}^{d}(\Delta),
       \qquad
       \sum_{j=1}^d I(\mu_j) \leq a_1 J_Q(\mu) + a_0 .
   \end{equation}
\end{proposition}
\begin{proof}
Consider a vector $y$ in the range of $C=B^tB$ that satisfies
\eqref{assumpH2}, and note that, since for all indices $i$,
$y_i^2>0$, the minimum $m=\min(y_i^2)$ is positive. Let
$x$ be a non-zero vector in $\R^r$ such that
$y=B^{t}x$, and $Q$ be an orthogonal matrix with $x/\|x\|$ as
its
first column. Then, the first row of $Q^{t}B$ is $y^{t}/\|x\|$
and
\[
J (\mu)=J (C\mu,\mu)= J
(Q^{t}B\mu,Q^{t}B\mu)=\frac{1}{\|x\|^{2}}I
\left(\sum_{j=1}^{d}y_{j}\mu_{j}\right)+
\sum_{k=2}^{r}I (\lambda_{k}),
\]
where we have set $(\lambda_{1},\ldots
,\lambda_{r})^{t}=Q^{t}B\mu$.
Next, we have the following lower bounds for the energies,
\begin{align}\notag
I (\mu_{j},\mu_{k}) & \geq
-\|\mu_k\|\int\log(1+|t|)d\mu_j(t)-\|\mu_j\|\int\log(1+|z|)d\mu_k(z)\\\label{ineq2}
& \geq -\gamma_{j,k}\left(\int Q_jd\mu_j+\int Q_kd\mu_k\right)
-\Gamma_{j,k},
\end{align}
where, as in the proof of Proposition \ref{lem_ineqQ-I}, the
positive real number $\gamma_{j,k}$ can be arbitrarily small
and $\Gamma_{j,k}$ is a sufficiently large number.
Hence, from \eqref{ineq2} applied with $j=k$, we deduce
$$
m\sum_{j}^{}\left(I (\mu_{j})+2\gamma_{j,j}\int
Q_jd\mu_j+\Gamma_{j,j}\right) \leq \sum_{j}^{}y_{j}^{2}
\left(I (\mu_{j})+ 2\gamma_{j,j}\int
Q_jd\mu_j+\Gamma_{j,j}\right),
$$
so that
\begin{align}\notag
m\sum_{j}I (\mu_{j})& \leq \sum_{j}^{}y_{j}^{2}I (\mu_{j})+
\sum_j(y_j^2-m)\left(2\gamma_{j,j}\int
Q_jd\mu_j+\Gamma_{j,j}\right)
\\
\notag
& = I \left(\sum_{j}^{}y_{j}\mu_{j}\right)-\sum_{j\neq
k}^{}y_{j}y_{k}I
(\mu_{j},\mu_{k})+
\sum_j(y_j^2-m)\left(2\gamma_{j,j}\int
Q_jd\mu_j+\Gamma_{j,j}\right)\\\notag
& =\|x\|^{2}J_Q (\mu)-2\|x\|^2\sum_{j=1}^d\int Q_jd\mu_j
-\|x\|^{2}\sum_{k=2}^{r}I (\lambda_{k})-\sum_{j\neq
k}^{}y_{j}y_{k}I
(\mu_{j},\mu_{k})\\
& +
\sum_j(y_j^2-m)\left(2\gamma_{j,j}\int
Q_jd\mu_j+\Gamma_{j,j}\right).
\label{ineq1}
\end{align}
For the signed measures $\lambda_k$, we have lower bounds
provided by Proposition \ref{lem_ineqQ-I},
\begin{equation}
I(\lambda_k)\geq -\gamma_k\sum_{j=1}^d\int Q_jd\mu_j
+\Gamma_k,\qquad k=1,\ldots,r.
\label{ineq3}
\end{equation}

Finally, for indices $j,k$ such that
$y_{j}y_{k}< 0$, we know from \eqref{assumpH2} that $\dist
(\Delta_{j},\Delta_{k})$ is positive so that in this case we
also have the upper bound,
\begin{equation}\label{ineq4}
I (\mu_{j},\mu_{k})  \leq \log
\left(\frac{1}{\dist (\Delta_{j},\Delta_{k})}
\right)\|\mu_{j}\|\|\mu_{k}\|.
\end{equation}
Making use of
\eqref{ineq2}, \eqref{ineq3} and \eqref{ineq4} into
\eqref{ineq1} leads to
$$
m\sum_{j=1}^{d}I (\mu_{j}) \leq
\|x\|^{2}J_Q (\mu)-c\sum_{j=1}^d\int Q_jd\mu_j-\Gamma,
$$
where $c$ is a positive real number since the constants
$\gamma_{k}$, $\gamma_{j,j}$, and $\gamma_{j,k}$ are
arbitrarily small. As the external fields $Q_j$ are lower
bounded, the above inequality implies
\eqref{ineqJI} with two constants $a_0$ and $a_1$ that depend
only on the tuple of sets $\Delta$, the interaction matrix $C$
and the compact set of masses $K$.
\end{proof}
\section{Existence of a solution}
\label{Existence}
In this section, we give the proof of Theorem \ref{existence}.
 We also prove, under an additional
technical assumption, that the components of a solution have compact
supports.
\begin{proof}[\bf Proof of Theorem \ref{existence}]
We show that $J_Q^*<+\infty$ as in \cite[Theorem I.1.3(a)]{SaTo}.
For
  $\epsilon>0$, the sets
$\Delta_j(\epsilon)=\{ x \in \Delta_j : Q_j(x) \leq 1/\epsilon
\}$
  are closed and thus compact by assumption on $Q_j$. Since $Q_j$ is finite on a set of positive
capacity,
$\Delta_j(\epsilon)$ is
of positive capacity for sufficiently small
$\epsilon>0$. Denoting by $\omega_{\Delta_j(\epsilon)}$ the
equilibrium
measure of such a
$\Delta_j(\epsilon)$ of positive capacity,
$I(\omega_{\Delta_j(\epsilon)})<\infty$, we find for the
$d$-tuple of measures $\mu\in\MM_K^d(\Delta)$ with
$\mu_j=b_j
\omega_{\Delta_j(\epsilon)}$, $j=1,...,d$, and $b=(b_j)\in
K$ that $J_Q(\mu)<\infty$.
Next, we prove that
$J^*_{Q}>-\infty$.
We have
\begin{equation*}
J_Q(\mu)=\sum_{k=1}^rI(\lambda_k)+2\sum_{j=1}^d\int Q_jd\mu_j,
\end{equation*}
where $(\lambda_1,\ldots,\lambda_r)^t=B\mu$, and the energies
$I(\lambda_k)$ satisfy inequalities of the type
(\ref{ineqQ-I}) with arbitrarily small positive constants
$\gamma_k$, the sum
of which can be made less than 1. Hence, there exists a constant
$\Gamma$ such that
\begin{equation}
\label{min-J_Q}
\forall \mu\in\MM_K^d(\Delta),\qquad
J_Q(\mu)\geq\sum_{j=1}^d\int Q_jd\mu_j-\Gamma\geq
-\sum_{j=1}^d |q_j|M_j(K)-\Gamma,
\end{equation}
with
$$q_j=\inf_{z\in\C} Q_j(z)>-\infty,\quad
M_j(K)=\sup_{\mu\in\MM_K^d(\Delta)}\|\mu_j\|<\infty,\qquad
j=1,\ldots,d.
$$
This finishes the proof of assertion (a).

The proof of assertion (b) follows
usual lines, see e.g. \cite[Chapter 5]{niso}, by constructing
$\mu^*$
as a weak limit of a minimizing sequence of $J_Q$. We first note,
in view of (\ref{min-J_Q}), that
for minimizing the energy $J_Q$, it is sufficient to consider
the subset $\TT$ of $\MM_K^d(\Delta)$ consisting of $d$-tuples of measures
$\mu$ such that
\begin{equation}\label{bound_Qmu}
\sum_{j=1}^d\int Q_jd\mu_j\leq J_Q^*+\Gamma+1.
\end{equation}
Let us show that $\TT$ is a log-tight family. For $\mu\in\TT$,
we have
$$
\sum_{j=1}^d\int (Q_j-q_j)d\mu_j\leq
J_Q^*+\Gamma+1+\sum_{j=1}^d|q_j|
M_j(K).
$$
We simply denote by $\alpha$ the right-hand side of the above
inequality. Let $\epsilon>0$ be given. Since the $Q_j$ are
admissible, there exists a compact set $K\subset\C$ such that
$$\sum_j(Q_j(x)-q_j)\geq\frac{\alpha}{\epsilon}\log(1+|x|),\qquad
x\in\C\setminus K.
$$
Consequently, for any $d$-tuple of measures $\mu$ in $\TT$,
$$\sum_j\int_{\C\setminus
K}\log(1+|x|)d\mu_j\leq\frac{\epsilon}{\alpha}\sum_j\int_{\C\setminus
K}(Q_j(x)-q_j)d\mu_j\leq
\frac{\epsilon}
{\alpha}\sum_j\int_{\C}(Q_j(x)-q_j)d\mu_j\leq\epsilon,
$$
which shows that the set $\TT$ is indeed log-tight.
Now, consider a minimizing sequence of $d$-tuples of 
  measures  $\mu^{(n)}\in \TT$, namely
  $$
         \lim_{n \to \infty} J_Q(\mu^{(n)}) = J^*_{Q}.
  $$
The family $\TT$ being log-tight, it is a fortiori tight, so
that by Theorem \ref{Prohorov}, there exists a subsequence,
that we still denote by $\mu^{(n)}$,
having a weak limit $\mu^*$. Its components $\mu_j^*$
are supported on $\Delta_j$, and its $d$-tuple of masses
belongs to $K$. Since $\log(1+|x|)$ is a continuous and lower bounded
function, we get from Lemma \ref{weak-lower-sc} that
$$\int\log(1+|x|)d\mu_j^*\leq
\liminf_{n\to\infty}\int\log(1+|x|)d\mu_j^{(n)},\qquad
j=1,\ldots,d.
$$
Moreover, up to an additive constant, $\log(1+|x|)$ is upper bounded by
$Q_j(x)$, inequality (\ref{bound_Qmu}) holds true for the
sequence $\mu_j^{(n)}$, and
$$
-|q_j|M_j(K)\leq q_j\|\mu_j^{(n)}\|\leq\int
Q_jd\mu_j^{(n)},\qquad
j=1,\ldots,d.
$$
Therefore, we may deduce that
$$\int\log(1+|x|)d\mu_j^*(x)<\infty,\qquad j=1,\ldots,d.
$$
Next, we show that each component $\mu_{j}^*$ is of finite
energy.
From Lemma \ref{scalar_lem_sc} follows that
$$
I(\mu_k^*)\leq\liminf_{n\to\infty} I(\mu_k^{(n)}),\qquad
k=1,\ldots,d.
$$
Adding these inequalities over $k$, and noting that, in view
of Proposition \ref{ineqJ-I}, the sum obtained on the right-hand
side is finite, we get
\begin{align}
I (\mu_{j}^*)
& \leq \liminf_{n\to\infty}\sum_{k=1}^{d}I (\mu_{k}^{(n)})-
\sum_{k\neq j}^{}I (\mu_{k}^*)
 \leq a_{1}\liminf_{n\to\infty} J_Q(\mu^{(n)})+a_{0}-
\sum_{k\neq j}^{}I (\mu_{k}^*) \notag\\
& = a_{1}J^{*}_{Q}+a_{0}-\sum_{k\neq
j}^{}I (\mu_{k}^*)<\infty, 
\end{align}
where the last inequality comes from
$$
I(\mu_k^*)\geq -2\|\mu_k^*\|\int\log(1+|x|)d\mu_k^*(x)
>-\infty\qquad k=1,\ldots,d.
$$

Consequently, $\mu^*\in \mathcal M_{K}^{d}(\Delta)$. From the
  lower semi--continuity of $J_Q$ established in
Proposition~\ref{lem_semicontinuity}, we conclude that
$J^*_Q\geq
J_Q(\mu^*)$, and thus $J^*_Q= J_Q(\mu^*)$, showing that
$\mu^*$ is
a minimizer of the extremal problem \eqref{min-pb}.
\end{proof}

  We now turn to the question of whether the supports of the
components of an extremal tuple of measures as in
Theorem~\ref{existence} are compact sets. This property was shown to hold
true under
  more restrictive conditions on the matrix $C$ and the tuple of sets $\Delta$ in \cite{BB,SaTo}.
In our generalized setting, we have the following result.

\begin{theorem}\label{sol-compact_bis}
Let $\mu\in\MM_K^d(\Delta)$ be a solution to the minimization
problem (\ref{min-pb}). Then, the components $\mu_i$,
$i=1,\ldots,d$, of $\mu$, have compact supports if and only if the following assertion holds true:
\\[\baselineskip]
there exists a real $\alpha$ and a number $M>0$ such that, for all pair
$(i,j)$ with $\Delta_i$ and $\Delta_j$ unbounded and
$c_{i,j}<0$,
there
holds
\begin{equation} \label{growth}
        c_{i,j} U^{\mu_j}(z) + \frac{1}{d} Q_i(z) \geq \alpha,
      \quad
\text{$\mu_i$--almost everywhere on $\Delta_i \setminus
D_M$,}\end{equation}
where $D_M$ denotes the
   closed disk of radius $M$ centered at zero.
\end{theorem}

\begin{remark}
The assumption (\ref{growth}) bears some similarity with
assumption [A3] in \cite[Definition 2.1]{BB}, where it is
assumed
that the functions
$$c_{i,j}\log\frac{1}{|z-t|}+\frac{Q_i(z)+Q_j(t)}{d},\qquad i,j=1,\ldots,d$$ are
uniformly lower bounded on $\Delta_i\times\Delta_j$.
\end{remark}
\begin{remark}\label{remark_NSatinf}
Assumption (\ref{growth}) is trivially satisfied if
\begin{equation*}
\forall i,j,\quad\text{if $\Delta_i$ and $\Delta_j$ are
unbounded
then $c_{i,j}\geq 0$}.
\end{equation*}
This condition can be seen as an analog of (\ref{compatNS})
where
we only consider the point at infinity in the intersection of
$\Delta_i$ and $\Delta_j$ (in the Riemann sphere). Of course, it is more
restrictive than the condition (\ref{growth}) but it has the
advantage that, for a given extremal problem, it can be
checked at
once from the data if it holds true or not.
\end{remark}
\begin{remark}
 Condition (\ref{growth}) follows from
   \eqref{assumpH2} in the case of a matrix $C$ of rank 1, for instance when considering a
   condenser as in \cite[Chapter~VIII]{SaTo}.
Indeed, here necessarily $C$ is a positive multiple of $y y^t$
with
the vector $y$ as in \eqref{assumpH2}. Thus $c_{i,j}<0$ implies
that
   $y_iy_j<0$, and hence for all $z \in \Delta_i$
   $$
         U^{\mu_j}(z) \leq \|
\mu_j\| \, \log
\left(\frac{1}{\mbox{dist}(\Delta_i,\Delta_j)}\right) .
   $$
   Consequently, \eqref{growth} follows by recalling
   that $Q_i$ is
   lower bounded. Hence,
   as in \cite[Theorem~VIII.1.4]{SaTo}, we may conclude that
   the components of an extremal tuple of measures in \eqref{min-pb}
   in the case $\text{rank}(C)=1$ have compact support.
\end{remark}
\begin{proof}
   Suppose first that the support of the measures
   $\mu_i$ are compact. Then, for $M$ sufficiently large, the sets $\supp(\mu_{i})\setminus D_{M}$
are empty sets so that (\ref{growth}) is trivially true.

   Conversely, let us show that $\mu_i$ has a compact support if
(\ref{growth}) holds.
   We first establish a property of $\mu_{i}$ similar to
   \cite[Lemma~5.4.1]{niso}, namely,
   \begin{equation} \label{stationnary}
\forall \nu_i \in \mathcal M_{\|\mu_i\|}(\Delta_i): \quad \int
(
U_i^{\mu} + Q_i ) d(\nu_i-\mu_i) \geq 0 .
   \end{equation}
For a proof of (\ref{stationnary}), we define $\nu\in
\mathcal M_K^d(\Delta)$ by $\nu_j=\mu_j$ for $j\neq i$. Notice
that $\mu + t(\nu-\mu) \in \mathcal
M_K^d(\Delta)$ for any $0<t<1$, and hence by definition of $\mu$
\begin{eqnarray*}
      0 &\leq& J_Q(\mu + t(\nu-\mu)) - J_Q(\mu)
      \\&=& 2 t \int ( U_i^{\mu} + Q_i ) d(\nu_i-\mu_i) + t^2
      J(\nu-\mu) .
   \end{eqnarray*}
   Dividing by $t$ and letting $t\to 0$ gives the desired
   inequality (\ref{stationnary}).

For our proof of compactness of $\supp(\mu_i)$, we may suppose
without
   loss of generality that $\Delta_i$ is unbounded, $\|
\mu_i \|>0$, and that $\mu_i(D_M)>0$, where for the last
property we
   possibly choose a larger $M$.
   We consider
   $$
\nu_i := \frac{\|\mu_i\|}{\mu_i(D_M)} \mu_i|_{D_M}   $$
being clearly an element of $\mathcal
M_{\|\mu_i\|}(\Delta_i)$. Then we may rewrite condition
(\ref{stationnary}) as
   $$
\Bigl( \frac{\|\mu_i\|}{\mu_i(D_M)} - 1 \Bigr)
\int_{|z|\leq M} ( U_i^{\mu} + Q_i ) d\mu_i
          -
          \int_{|z|> M}  ( U_i^{\mu} + Q_i ) d\mu_i \geq 0,
   $$ or
   \begin{equation} \label{compact}
        (\|\mu_i\| - \mu_i(D_M)) \alpha_0 \geq
        \|\mu_i\| \int_{|z|> M}
        ( U_i^{\mu}(z) + Q_i(z) ) d\mu_i(z) ,
   \end{equation}
   with the finite constant
   $$
        \alpha_0 := \int ( U_i^{\mu}(z) + Q_i(z) ) d\mu_i(z) .
   $$
It remains to show that $U_i^{\mu}(z) + Q_i(z)$ is
sufficiently large
   $\mu_i$--almost everywhere on $\Delta_i \setminus D_M$.

   For this, notice first
   that, by possibly making $\alpha$ smaller and $M$ larger,
(\ref{growth}) also holds for all indices $j$ with $c_{i,j}<0$
and compact $\Delta_j$ since then $\supp(\mu_j)$ is compact.
   In case $c_{i,j}\geq 0$ we use (\ref{pot-lower-bd}) to
   conclude that, for
   $\mu_i$--almost all
   $z\in \Delta_i\setminus D_M$,
   \begin{eqnarray*}
U_i^{\mu}(z) + Q_i(z) &\geq& \sum_{j,c_{i,j}\geq 0}
c_{i,j}U^{\mu_j}(z) + \sum_{j,c_{i,j}< 0} ( \alpha -
\frac{1}{d}
Q_i(z)) + Q_i(z)
       \\&\geq& - \sum_{j,c_{i,j}\geq 0} c_{i,j} \| \mu_j \|
       \log(1+|z|) + \frac{1}{d} Q_i(z) + \alpha_1
   \end{eqnarray*}
for some constant $\alpha_1$. Here we have used the fact
that $c_{i,i}\geq 0$. According to the third condition of
admissibility
in
Definition \ref{def-adm}, i.e. the behavior of $Q_i$ at infinity, we
may now possibly choose a larger $M$ such that
   $U_i^{\mu}(z) + Q_i(z) \geq (\alpha_0+1)/\|\mu_i\|$ for $\mu_i$--almost all
   $z\in \Delta_i\setminus D_M$.
   Hence inequality
   (\ref{compact}) becomes
   $$
        (\|\mu_i\| - \mu_i(D_M)) \alpha_0 \geq
        (\|\mu_i\| - \mu_i(D_M)) (\alpha_0+1),
   $$
   implying that $\|\mu_i\| = \mu_i(D_M)$, and the fact that
   $\mu_i$ has compact support.
\end{proof}
\section{Uniqueness and equilibrium conditions}\label{Unique}

\begin{proof}[{\bf Proof of Theorem~\ref{uniqueness}(a)}]
Our proof relies on
Proposition~\ref{lem_bounded}, but otherwise the arguments of
\cite{niso} or \cite[Proof of Theorem~1.1]{AL} carry over to
our more general setting.

 It is sufficient to show that
  the application $\mu \mapsto J_Q(\mu)$ is strictly
convex\footnote{More precisely, we only show strict midpoint
convexity, which is sufficient for our purposes. However,
   together with
the lower semi--continuity established in
Proposition~\ref{lem_semicontinuity} one
   may deduce strict convexity.}
in the convex subset of $\MM^d_K(\Delta)$ consisting of $d$-tuples of 
measures $\mu$ with finite $J_Q$--energy. By finiteness of the
$J$--energy on $\MM^d_K(\Delta)$, that simply boils down to
    $$\int Q_jd\mu_j <\infty,\qquad j=1,\ldots,d.
    $$
For two distinct $d$-tuples of measures $\mu$ and $\nu$ of finite
$J_Q$-energies, we have
   \begin{eqnarray*} &&
         \frac{1}{2}\Bigl( J_Q(\mu) + J_Q(\nu)
         \Bigr) -
         J_Q\left(\frac{\mu + \nu}{2}\right)
         = \frac{1}{2}\Bigl( J(\mu) + J(\nu)
         \Bigr) -
         J\left(\frac{\mu + \nu}{2}\right)
        \\&& =
        \sum_{i,j=1}^d c_{i,j} \Bigl(
\frac{1}{2}\Bigl( I(\mu_i,\mu_j) + I(\nu_i,\nu_j) \Bigr) -
I\left(\frac{\mu_i + \nu_i}{2},\frac{\mu_j + \nu_j}{2}\right)
\Bigr)
         = J\left(\frac{\mu - \nu}{2}\right),
   \end{eqnarray*}
   and it only remains to show that the last term is positive.
By the definition \eqref{def-K} of $K$, the vector of masses
$\|\mu \| - \| \nu \|$ is an
   element of the kernel of the matrix $A$, which by
\eqref{def-C} and \eqref{image_AC} is a subset of the kernel of
$C$ and thus of $B$.
Hence the strict convexity follows from Proposition
\ref{lem_bounded}.
\end{proof}

\begin{remark}
  There are other sufficient conditions to ensure strict
  convexity of the map $\mu \mapsto J_Q(\mu)$ on $d$-tuples of measures of
$\MM^d_K(\Delta)$ of finite $J_Q$-energy, for instance we
may replace \eqref{image_AC} by the requirement that the union of
the
  $\Delta_j$ is compact, with capacity less than 1.
  Another sufficient condition for strict convexity, namely
  \begin{equation}\label{cij0}
    \forall i\neq j, \quad \mbox{if~~~}
\Delta_i\cap\Delta_j\neq\emptyset
\mbox{~~~then~~~}
c_{ij}=0,
  \end{equation}
has been considered in \cite{AL,GRS}. Notice that \eqref{cij0}
is
stronger than \eqref{compatNS}, and that
\eqref{compatNS} alone
  does not imply
  strict convexity, see Example~\ref{example0_not_uniqueness}.
\end{remark}

\begin{remark}\label{imCimAT}
   We claim that if there is equality in assumption
   \eqref{image_AC} then \eqref{assumpH2} holds. To see this,
notice that from the full rank decomposition $C=B^tB$ and from
the assumption
$\Im(C)=\Im(A^t)$ we conclude that there exists a matrix $E$ such
that $A=E B$, implying that we may rewrite the non empty compact
$K$ as $K=\{x\in\R^d_+, Bx= b \}$ for a suitable
vector $b \in \mathbb R^r$. Writing $e=(1,...,1)^t\in \mathbb
R^d$, we conclude that the linear optimization problem $\max \{
e^t x ,Bx=b, x \geq 0 \}$ has an optimal solution. In particular \cite[Theorem~19.12]{BGLS}, there is a Lagrange multiplier $\lambda \in \mathbb R^r$ with $B^t \lambda \geq e$. Hence $y:=B^t \lambda$
is an element
of $\Im(C)=\Im(B^t)$ with strictly positive components,
implying
\eqref{assumpH2}.
\end{remark}

   Before entering the details of the proof of assertion (b) of
Theorem~\ref{uniqueness},
  we shortly comment on the equilibrium
  conditions \eqref{char-solQC1} and \eqref{char-solQC2}. First recall from
Lemma~\ref{lem-tsuji}
that the potentials $U^{\mu_j}$ for $j=1,...,d$ are finite and hence
$U_i^\mu+Q_i$
is well--defined in $\Delta_{i} \setminus \Delta_{i,\infty}$
with
  $\Delta_{i,\infty}\subset \Delta_i$ some polar Borel set.
  Also, $U_i^\mu+Q_i$ as
  a sum of measurable functions is measurable, and hence both sets
\begin{align*}
\Delta_{i,+} & = \{ z\in \Delta_i \setminus \Delta_{i,\infty},~
U_{i}^{\mu} (z)+Q_{i} (z)> (A^t F)_i \},\\
\Delta_{i,-} & = \{ z\in \Delta_i \setminus \Delta_{i,\infty},~
U_{i}^{\mu} (z)+Q_{i} (z)< (A^t F)_i \},
\end{align*}
are Borel sets.
  Hence \eqref{char-solQC1} means that
  $\Delta_{i,\infty} \cup \Delta_{i,-}$ is polar,
  whereas \eqref{char-solQC2} can be equivalently rewritten as
  $\mu_i(\Delta_{i,\infty} \cup \Delta_{i,+})=0$.

As
in \cite[Lemma~5.4.2]{niso} we have to establish a
different characterization of an extremal tuple of measures which
generalizes
(\ref{stationnary}).
\begin{lemma}\label{convex_problem}
The $d$-tuple of measures $\mu=(\mu_1,...,\mu_d)\in
\MM^d_K(\Delta)$
with $J_Q(\mu)<\infty$ is extremal for \eqref{min-pb} if
and only
if for any $d$-tuple of measures $\nu=(\nu_1,...,\nu_d)\in
\MM^d_K(\Delta)$
    with $J_Q(\nu)<\infty$ we have
    \begin{equation} \label{variation1}
        \sum_{i=1}^d \int (U_i^\mu + Q_i) d\nu_i
        \geq
        \sum_{i=1}^d \int (U_i^\mu + Q_i) d\mu_i .
    \end{equation}
\end{lemma}
\begin{proof}
In order to see that \eqref{variation1} is necessary for
optimality, notice
   that, for all $0<t\leq 1$, we have $\mu + t (\nu-\mu)\in
   \MM^d_K(\Delta)$, with
   \begin{equation} \label{variation2}
                  J_Q(\mu + t (\nu-\mu))-J_Q(\mu)
                  = 2 t
                         \sum_{i=1}^d \int (U_i^\mu + Q_i)
                         d(\nu_i-\mu_i) + t^2 J(\nu-\mu)
   \end{equation}
   being nonnegative. Dividing by $t$ and letting $t\to 0$ gives
   \eqref{variation1}. Conversely, we recall
from Proposition~\ref{lem_bounded} that $J(\nu-\mu)\geq 0$.
Injecting
\eqref{variation1} into \eqref{variation2} for $t=1$, we
conclude as required
   that $\mu$ is extremal.
\end{proof}
\begin{proof}[{\bf Proof of Theorem~\ref{uniqueness}(b)}]
   Suppose first that $\mu\in \MM^d_K(\Delta)$ satisfies
   \eqref{char-solQC1} and \eqref{char-solQC2}.
   Then $\mu_i(\Delta_{i,\infty} \cup \Delta_{i,+})=0$,
   and integrating \eqref{char-solQC2}
   with respect to $\mu_i$ shows that $J_Q(\mu)<\infty$.
   Let now $\nu\in \MM^d_K(\Delta)$ with $J_Q(\nu)<\infty$.
   Then $\nu_i(\Delta_{i,\infty} \cup \Delta_{i,-})=0$ by
   Lemma~\ref{mes_zero}.
Hence integrating \eqref{char-solQC1} with respect to $\nu_i$
and \eqref{char-solQC2} with respect to $\mu_i$ gives
   $$
        \sum_{i=1}^d \int (U_i^\mu + Q_i) d\nu_i
        - \sum_{i=1}^d \int (U_i^\mu + Q_i) d\mu_i
\geq \sum_{i=1}^d (\| \nu_i \| - \| \mu_i \| ) (A^t F)_i = 0,
   $$
   the last equality following from the definition
   of the polyhedron of masses $K$.
Hence $\mu$ is extremal according to
Lemma~\ref{convex_problem}.

Suppose now that
$\mu\in\MM^d_K(\Delta)$ is extremal. Consider the set of indices $I = \{ i\in \{1,...,d\} : \| \mu_i\| >0 \}$, set for $i \in I$
   $$
         w_i := \frac{1}{\|\mu_i\|} \int (U_i^\mu+Q_i)d\mu_i ,
   $$
   and consider as before the Borel sets
$\Delta_{i,+}= \{ z\in \Delta_i \setminus \Delta_{i,\infty} :
U_{i}^{\mu} (z)+Q_{i} (z)> w_i \}$ and
$\Delta_{i,-}= \{ z\in \Delta_i \setminus \Delta_{i,\infty} :
U_{i}^{\mu} (z)+Q_{i} (z)< w_i \}$.
Following
\cite[Proposition~5.4.4]{niso}, we claim that, for
$i\in I$,
   \begin{equation} \label{equi1}
U_i^\mu(x)+Q_i(x) \geq w_i, \mbox{~~quasi--everywhere on $\Delta_i$}.
\end{equation}
Suppose the contrary for some $i\in I$. Since $\Delta_{i,\infty}$
is polar,
we conclude that $\Delta_{i,-}$ is of positive capacity. Thus
there
   exists a compact
   set $E\subset\Delta_i$ with $U_i^\mu$ well-defined
   and finite
   on $E$,
   $\capa(E)>0$,
   and $U_i^\mu(x)+Q_i(x) < w_i$ for all $x\in E$. Taking
   any $\nu_i\in\MM_{\| \mu_i
   \|}(E)$, then with $\nu_j=\mu_j$ for $j\neq i$ we get
$\nu\in \MM_K^d(\Delta)$ and, by Lemma~\ref{convex_problem}, $$
           0 \leq
\sum_{\ell=1}^d \int (U_\ell^\mu + Q_\ell) d(\nu_\ell-\mu_\ell)=
\int (U_i^\mu+Q_i)d\nu_i - \| \nu_i \| w_i ,
   $$
but the term on the right is negative by construction of $E$ and
$\nu_i$, a contradiction. Thus \eqref{equi1} holds.

   Following
\cite[Proposition~5.4.5]{niso}, we now claim that, for
$i\in I$,
   \begin{equation} \label{equi2}
U_i^\mu(x)+Q_i(x) \leq w_i, \mbox{~~$\mu_i$--almost
everywhere.}\end{equation}
   Suppose the contrary for some $i\in I$.
Since $\mu_i(\Delta_{i,\infty})=0$ by
Lemma~\ref{mes_zero},
we get $\mu_i(\Delta_{i,+})>0$. Applying, e.g.,
   \cite[Theorem~2.18]{rudin}, we conclude that
     there exists a compact
   set $E\subset\Delta_i$ with $U_i^\mu$ well-defined
   and finite on $E$, $\mu_i(E)>0$,
   and $U_i^\mu(x)+Q_i(x) > w_i$ for all $x\in E$.
A combination of Lemma~\ref{mes_zero} with \eqref{equi1} tells
us that
   $$
\int_{\Delta_i\setminus E} (U_i^\mu+Q_i) d\mu_i \geq w_i
\mu_i(\Delta_i\setminus E),
   $$
   and thus
   $$
         \| \mu_i \| w_i \geq \int_E (U_i^\mu+Q_i) d\mu_i
         + w_i \mu_i(\Delta_i\setminus E) >
                 w_i \mu_i(E)
         + w_i \mu_i(\Delta_i\setminus E) ,
   $$
a contradiction. Hence also (\ref{equi2}) is true. Thus we
have shown so far that,
for indices $i$ with $\|\mu_i\|>0$, \eqref{char-solQC1} and
\eqref{char-solQC2} hold true if
we replace $(A^tF)_i$ by a suitable constant $w_i\in \mathbb
R$. It remains to relate these constants $w_i$ with $A$ and also
to discuss the partial potentials $U^\mu_i$ for indices $i$
such that $\| \mu_i \|=0$.
For this purpose, similar to \cite[Part 3 of proof of
Theorem~1.2]{AL},
we consider the quadratic optimization problem in $\mathbb
R^d$,
$$
          \min \{ x^t H x + 2h^t x,\quad x \in K \},
   $$
    where $H\in\R^{d\times d}$ and $h\in\R^d$ with
    $$
    H_{i,j} = c_{i,j} I(\nu_i,\nu_j), \quad
          h_i = \int Q_i \, d \nu_i ,\qquad i,j=1,\ldots,d
                    $$
   and the probability measures  $\nu_i\in \MM_{1}(\Delta_i)$
   are defined by $\nu_i=\mu_i/\| \mu_i\|$ if
   $\|\mu_i\| \neq 0$, and else arbitrary but fixed. Then, by
   Theorem~\ref{uniqueness}(a), $\| \mu \| \in K$ is the
   unique solution of the above quadratic problem. From
   \cite[Theorem~19.12]{BGLS} we know
   that there exist Lagrange multipliers $F\in \mathbb R^m$ and
   $G\in \mathbb R^d$ such that
   \begin{equation}
    \label{KKT!}
      H \| \mu \| + h = A^t F + G , \quad
\forall i , \quad G_i \geq 0 , \quad \| \mu_i \| \, G_i = 0 .
\end{equation}
   In case $\| \mu_i \|\neq 0$ we find from
   \eqref{equi1} and \eqref{equi2} that
   \begin{eqnarray*} &&
         (H \| \mu \| + h)_i =
\sum_{j=1}^d c_{i,j} I (\frac{\mu_i}{\|\mu_i\|},\mu_j) + \int
Q_i \, \frac{d\mu_i}{\|\mu_i\|}
         = \int (U^\mu_i+Q_i) \frac{d\mu_i}{\|\mu_i\|} = w_i .
   \end{eqnarray*}
   Also, $G_i=0$, and hence
   $(H \| \mu \| +h)_i = w_i=(A^tF)_i$. In particular,
   relations    \eqref{equi1} and \eqref{equi2} imply
   the desired relations \eqref{char-solQC1} and
   \eqref{char-solQC2}.
   In case $\| \mu_i \| = 0$
   we learn from \eqref{KKT!} that
   $$
       \forall \nu_i \in \MM_{1}(\Delta_i) ,
       \quad
\int \Bigl( U_i^{\mu} + Q_i \Bigr) d\nu_i \geq (A^tF)_i,   $$
implying \eqref{char-solQC1}, and assertion
\eqref{char-solQC2} is trivially true.
\end{proof}
\section{A review of some examples}
\label{Examples}

Many recently studied problems, such as e.g. the behavior of
Hermite-Pad\'e
approximants, the limit eigenvalue distribution of banded
Toeplitz matrices, or the limit distribution of non-intersecting
brownian
paths, translate into vector equilibrium problems with
external fields.
Existence and uniqueness of the
solution were shown under conditions that are actually covered by the results of
the previous sections. The above-mentioned equilibrium problems can be
stated
in terms of graphs. We recall that for a graph  $G=({\cal V},
{\cal E})$, the set of edges  ${\cal E}$ is a subset of the cartesian product ${\cal
V}\times {\cal V}$, where ${\cal V}$ denotes the set of
vertices.
For multigraphs we allow for repeated edges between two given
vertices. We also remind the reader that the incidence matrix $A$ is
labelled in rows by vertices and in columns by edges, with a column
corresponding to an
edge from the vertex $u$ to the vertex $v$ has entry $-1$ at
row
$u$, $1$ at row $v$ and $0$ elsewhere.

\begin{center}
\begin{figure}[ht]
\centering{\includegraphics[scale=0.6]{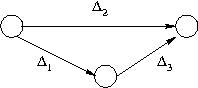}} \caption{A
graph
with undirected cycle but no directed cycle}\label{AV}
\end{figure}
\end{center}

In what follows, we always suppose that a graph or a
multigraph $G=({\cal V}, {\cal E})$ is given. We denote its
incidence
matrix
by $A$ and we consider as interaction matrix the matrix
$C=A^tA$, together with the polyhedron of masses
$K=\conj{x\in\R^d_+,~ Ax=a}$. In what follows $K$ is supposed to contain at least one element with strictly positive components. For instance, for the graph
of figure \ref{AV}, we have
\begin{equation}\label{exemple_AV}
 A=\left( \begin{array}{ccc}
-1 & -1 & 0\\ 0 & 1 & 1\\ 1 & 0 & -1
\end{array}
 \right),\quad C=\left( \begin{array}{ccc}
2 & 1 & -1\\ 1 & 2 & 1\\ -1 & 1 & 2\end{array}
 \right).\end{equation}
As a consequence, the interaction matrix $C$ is indexed in rows and columns
by the edges and it can be checked that its entries are
$-2,-1,0,1,2$ with the following interpretation
\begin{equation*}
C_{\alpha , \beta}=\left\{
\begin{array}{cl}
2& \mbox{ if }   \alpha=\beta \mbox{ or }~
\raisebox{-0.4cm}{\includegraphics[scale=0.6]{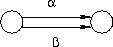}}
\\
1& \mbox{ if } ~ \includegraphics[scale=0.6]{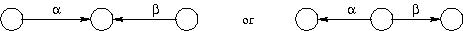}
\\
-1 & \mbox{ if } ~ \includegraphics[scale=0.6]{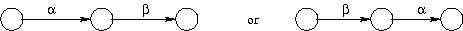}
\\
-2 &  \mbox{ if } ~ 
\raisebox{-0.4cm}{\includegraphics[scale=0.6]{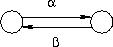}}
\\
0& \mbox{ elsewhere (i.e.~} \alpha, \beta \mbox{ do not have
any vertex in common)}
\end{array}\right.
\end{equation*}
By construction, the matrix $C$ is always positive
semi-definite.
To each edge $i$ we associate a closed set $\Delta_i$ and a
measure $\mu_i$ supported on $\Delta_i$.

We can interpret the different assumptions we made in the previous sections about the matrix $C$
and
the
supports $\Delta_i$ in terms of graph
theory.

\begin{proposition}\label{graphprop}
The following assertions hold true:
\begin{enumerate}
\item[(a)] The following three statements are equivalent: (i) matrix $C$ is invertible;
    (ii) $G$ has no undirected cycle;
     (iii) the polyhedron of masses $K =\conj{x\in\R_+^d;~ Ax=a}$ 
     is a singleton.
\item[(b)] The polyhedron of masses $K$ is compact if and only if $G$ has no directed cycle.
\item[(c)]
Condition (\ref{compatNS})
is equivalent to the fact that any two edges which follow each other
correspond to non intersecting sets
 $\Delta_i$ and $\Delta_j$.
\item[(d)]
Condition (\ref{cij0})
is equivalent to the fact that any two distinct edges
corresponding to intersecting sets $\Delta_i$ do not have any
vertex in common.
\item[(e)] Condition (\ref{assumpH1})
is equivalent to: $$\forall\mbox{ set $I$ of edges of }{\cal
E}\mbox{ forming an undirected cycle in } G, ~ \capa
(\cap_{\alpha \in I}\Delta_{\alpha} )=0.$$
\item[(f)] Let $G^*$ be the undirected intersection graph of
the sets $\conj{\Delta_i}_{i=1}^d$ that is, the vertices of $G^*$ are the
edges of $G$ and there is an edge in $G^*$ between $i$ and $j$
if the corresponding sets $\Delta_i$ and $\Delta_j$ are
intersecting. Condition (\ref{assumpH2})
is equivalent to:

each connected component of $G^*$ corresponds to a subgraph in
$G$
that does not contain  a directed cycle.

\end{enumerate}
\end{proposition}
We do not present here complete proofs for these assertions
which
follow from graph theory. Notice however that (a) is based on the classical fact that the rank of an incidence matrix is given by the number of its columns iff the underlying graph has no undirected cycle. Assertions
(c) and (d)
immediately follow from the above graph interpretation of the
entries of $C$.

Condition (\ref{cij0}) 
is obviously stronger than (\ref{compatNS}).
From the graph theory
interpretation given in assertions (d) and (e)
we
see that (\ref{cij0}) implies (\ref{assumpH1}). From
assertions (b) and (f) we see that the
compactness of the polyhedron $K$ implies (\ref{assumpH2}), as
noticed already in
Remark \ref{imCimAT}.

\begin{figure}[h]
\centering{\includegraphics[scale=0.6]{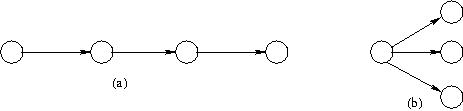}}
\caption{Tree graphs}\label{niki}
\end{figure}

The first vector equilibrium problems using the terminology of graphs were studied
in
\cite{GRS}, where systems of Markov
functions
 generated by a rooted tree   $G=({\cal V},\cal{E})$,
 the so-called {\em generalized Nikishin systems}, were considered.
Recall that a tree is a connected graph without undirected
cycles.
 In particular, by properties (a), (e), and (f) of
 Proposition \ref{graphprop}, $C$  is invertible and
conditions (\ref{assumpH2}) and (\ref{assumpH1}) are
satisfied. So the result \cite[Theorem 1]{GRS} also follows from
our work,
and we may drop in \cite[Theorem 1]{GRS} any further
requirements
 on the sets $\Delta_j$ like (\ref{cij0}) or (\ref{compatNS}).
 The authors associate to each vertex in ${\cal V}$ a
 Markov function, and to each edge $\alpha$ in ${\cal E}$
 a measure with support in an interval $\Delta_{\alpha}$. This
class includes the well-known Nikishin systems, see Figure
\ref{niki} (a),
and the Angelesco systems, see Figure \ref{niki} (b), with interaction matrices $C$ respectively
given by
$$
\begin{pmatrix}
2 & -1 & 0\\
-1 & 2 & -1\\
0 & -1 & 2
\end{pmatrix},\quad
\begin{pmatrix}
2 & 1 & 1\\
1 & 2 & 1\\
1 & 1 & 2
\end{pmatrix}.
$$
The solution of their extremal problem is related to the limit
distributions
of the zeros of the polynomial denominators of
the Hermite-Pad\'e approximants to the generalized Nikishin
systems.

\begin{figure}[ht]
\centering{\includegraphics[scale=0.6]{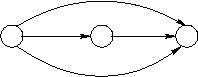}}
\caption{Rooted multigraph with undirected cycle}\label{ApteLysov}
\end{figure}

In  \cite{AL}, the results of \cite{GRS} were generalized to
rooted multigraphs $G=({\cal V},\cal{E},{\cal O}) $ with a
root ${\mathcal O}$, that is, multigraphs which have no
directed cycles
but do have directed  paths from ${\cal O}$ to any other vertex.  An
example of such a graph with undirected cycles is shown in
Figure
\ref{ApteLysov}.
By generalizing the
ideas of \cite{GRS}, the graph is associated to a system of Markov
functions with intersecting supports. According to assertion (f) of Proposition \ref{graphprop},
 condition (\ref{assumpH2}) holds  since there are no directed
cycles. Also, as said before, the condition (\ref{cij0})
imposed in \cite{AL} implies (\ref{assumpH1}). Thus
\cite[Theorem 1.1]{AL}
dealing with $K$ as in (\ref{def-K}) is covered by our work as
well. The Hermite-Pad\'e approximants to specific systems of
Markov functions related to graphs with cycles were
also investigated
in \cite{So} in connection with applications to number theory.

Another vector equilibrium problem appears in \cite{A} and
\cite{AKWa} in the study of the asymptotics of diagonal
simultaneous Hermite-Pad\'e approximants to two analytic functions with
separated pairs of branch points. The authors define the class
${\cal H} (\C\setminus \Gamma )$ of holomorphic functions in
$\C\setminus \Gamma$, where $\Gamma$ is a piecewise analytic
arc joining two points $a$ and $b$ in $\C$. A typical example of such
a function is
$$f(z)=\log\left(\frac{z-a}{z-b}\right).$$
For
$f_1\in
H(\C\setminus\Gamma_1)$,  $f_2\in H(\C\setminus\Gamma_2)$, with
 $$
 \Delta_1=\Gamma_1, \quad\Delta_{2}=
\mbox{Clos}(\Gamma_2\backslash\Gamma_1),
$$
and $\Delta_3$ a piecewise analytic arc containing the intersection $\Delta_1\cap \Delta_2$, they
show the existence and uniqueness of a triple of
measures
$$\mu=(\mu_1,\mu_2,\mu_3) \text{ with }\supp(\mu_i)\subset
\Delta_i,\quad i=1,2,3,$$
minimizing the energy $J(\mu)$, where
the interaction
matrix $C$ is given in (\ref{exemple_AV}), corresponding to the
graph in Figure \ref{AV}, and the set of masses is given by
$$
  K=\conj{x\in\R^3_+,~
  \left(\begin{array}{ccc}1& 1& 0\\ 1 & 0 &
-1\end{array}\right)\left(\begin{array}{c}x_1\\x_2\\x_3\end{array}\right)=\left(\begin{array}{c}
  2\\1\end{array}\right)}=\conj{x\in\R^3_+,~A x=
  \left(\begin{array}{c} -2\\ 1\\ 1\end{array}\right)}.
$$
Notice that this graph contains an undirected cycle but, since
$\capa (\Delta_1\cap \Delta_2)=0$, we are again in the settings
of
our theorems. The measure  $\mu_1+\mu_2$ is the  limit
distribution of  the poles of the diagonal simultaneous
Hermite-Pad\'e approximants of the functions $(f_1,f_2)$, and
the measure $\mu_3$ describes the limit distribution of the extra
interpolation points to $f_1$.

In \cite{DK}, the limit distribution of non-intersecting
one-dimensional Brownian paths with prescribed starting and
ending points is also characterized by a vector equilibrium
problem. As explained in \cite{DK}, there is an underlying
undirected graph $G_u$ whose edges connect vertices in
the
set of starting points with vertices in the set of ending
points,
that is, a bipartite graph. The authors show,
in addition, that their graph is a tree, see \cite[Proposition
2.1]{DK}. In
\cite[Corollary 2.9.]{DK}, they establish existence and
uniqueness
of a solution to an extremal vector equilibrium problem with
interaction matrix $C = (B^t B)/2$, $B$ being the incidence
matrix of $G_u$, with quadratic external fields, fixed masses,
and sets
$\Delta_j=\mathbb R$. The supports of the
extremal measures are compact, and
describe the limiting behavior of such
non-intersecting one-dimensional Brownian paths.
\begin{figure}[ht]
\centering{\includegraphics[scale=0.6]{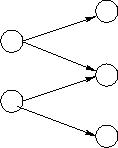}}
\caption{bipartite directed graph}\label{bipart}
\end{figure}
In order to relate  \cite[Corollary 2.9.]{DK} to our findings,
notice that, by considering the natural orientation of edges
from
starting points to ending points, we get a directed
graph $G$ which is both a tree and a bipartite graph, see the
example in
Figure \ref{bipart}. Using this last property, it is not
difficult
to see that $B^t B=A^t A$ with $A$ the incidence matrix of
$G$. Thus, we learn from assertion (a) of Proposition \ref{graphprop}
that $C$ is invertible, see also \cite[Proposition 2.8]{AL},
and
that $K$ is a  singleton. In particular, both conditions
(\ref{assumpH2}) and (\ref{assumpH1}) are true  and even the
condition (\ref{compatNS}) holds. Nevertheless, \cite[Corollary
2.9.]{DK} is not a consequence of \cite[Chapter 5]{niso} since the sets
$\Delta_j$ are not compact. However, the quadratic external
fields
of \cite{DK} are admissible in the sense of our Definition
\ref{def-adm}, and thus existence,
uniqueness and equilibrium conditions for an extremal tuple of measures
also
follow from our general findings. Note also that the
compactness of the supports of these extremal measures follows
from
Remark~\ref{remark_NSatinf} since all entries of $C$ are non
negative.

\section{Conclusion}\label{Conclusion}

    In this paper we have shown existence and uniqueness of an
    extremal tuple of measures for a vector generalization of
    a weighted
    energy problem in logarithmic potential theory
    with a polyhedron
of masses, substantially weakening the assumptions
typically assumed in other papers on this subject. We have also
derived a characterisation of such an extremal tuple of measure in
terms of
    equilibrium conditions for the vector potentials.

We have not been able to prove in our general setting that the
supports
of the components of the extremal tuple of measure are
always compact. We conjecture that, because of the growth of the external field at infinity and condition (\ref{assumpH2}), it should be true. In any case, we note that the variational inequality (\ref{char-solQC2}) implies that the potentials
$U^{\mu_{j}}$ such that $c_{i,j}<0$ satisfy
$
U^{\mu_{j}}(z)/\log|z|\to\infty$ as $z\in\Delta_{i}$ tends to infinity (up to a set of $\mu_{i}$-measure
zero). Hence, in view of assertion (ii) of \cite[Theorem 5.7.15]{Hel}, we may at least conclude that the support of $\mu_{i}$ is the union of a set of $\mu_{i}$-measure
zero and a set thin at infinity.


   There are also examples of vector--valued extremal problems
   in logarithmic potential theory where the external fields have a slow increase near $\infty$, or are
even not present. For instance, in \cite{duits_kuijlaars}, the
authors describe the limiting eigenvalue distribution of banded
Toeplitz matrices. It is obtained as a component
   of the solution of a vector equilibrium problem
	with a positive definite interaction matrix $C$ (namely the one
   of a Nikishin system),
   without any external field at all. Also,
in \cite{DelDui}, these results have been extended to Toeplitz
matrices with rational symbol, and in this case the vector
   equilibrium problem includes external fields of the form
   $Q(z)=C\log(|z|)$. In these examples, it may happen that
   the extremal measures do not have a bounded support.
For a general analysis of such examples, one should work on
the Riemann
   sphere instead of the complex plane, see the recent contribution \cite{arno_adrien} in case of positive definite $C$.
\\[\baselineskip]
\noindent{\bf Acknowledgements.} The authors would like to thank the reviewers for their careful reading of the manuscript. V. Kalyagin acknowledges supports by RFFI 10-01-00682 and Scientific School 8033.2010.1.

\appendix

\section{Appendix}\label{Proofs}
\begin{proof}[Proof of Lemma \ref{mes_zero}]
Assume $E$ is a Borel set such that $\mu(E)>0$. By regularity
of $\mu$, there exists a compact subset $K$ of $E$ with
$\mu(K)>0$. Set $\widetilde\mu=\mu_{|K}$. Then,
\begin{eqnarray*}
I(\widetilde\mu)& = & I(\mu)+\int_{\C\setminus
K}\int_\C\log(|z-t|)d\mu(z)d\mu(t) +\int_K\int_{\C\setminus
K}\log(|z-t|)d\mu(z)d\mu(t)\\ & \leq & I(\mu)
+4\int_\C\int_\C\log(1+|t|)d\mu(t)d\mu(z),
\end{eqnarray*}
which shows that $I(\widetilde\mu)<\infty$ and thus
$\capa(E)>0$.
\end{proof}
    \begin{proof}[Proof of Lemma \ref{weak-lower-sc}]
By \cite[Theorem 2.1.3]{ran}, there exists an increasing
sequence of continuous functions $h_m$ which converges
pointwise to $Q$. Assume $Q$ is lower bounded by $c\in\R$.
Set
    $$
    {\widetilde h}_m=\min(c+m,\max(c,h_m)).
    $$
Then, $({\widetilde h}_m)_m$ is an increasing sequence of
continuous bounded functions that still tends pointwise to $Q$
and we have
    $$
    \liminf_{n\to\infty}\int Qd\mu_n  \geq\lim_{m\to\infty}
    \liminf_{n\to\infty}\int {\widetilde h}_md\mu_n
     =\lim_{m\to\infty}\int {\widetilde h}_md\mu=\int Qd\mu,
    $$
where in the last equality we use the monotone convergence
theorem.
    \end{proof}
\begin{proof}[Proof of Lemma \ref{scalar_lem_sc}]
    Let $\epsilon>0$ be given and let $M>1$ be such that
$$\forall n\geq 0,\quad\iint_{|x-y|\geq M}\log(|x-y|)d\mu_n(x)
d\nu_n(y)\leq \epsilon. $$ Note that the existence of $M$
follows
from the simple inequalities
$$0\leq\log(|x-y|)\leq\log(1+|x|)+\log(1+|y|),$$ satisfied for
$|x-y|\geq1$, the fact that the masses of the measures are
uniformly bounded, and the log-tightness of the sequences. We
also
set $h(t)$ for a continuous function on $\R_+$ such that $$
0\leq
h(t)\leq 1,~\forall t\in\R_+,\quad h(t)=1\text{ for }t\leq
M,\quad
h(t)=0\text{ for }t\geq M+1. $$ Then, we have
\begin{align}\notag
I(\mu_n,\nu_n)& =\iint
\log(|x-y|^{-1})h(|x-y|)d\mu_n(x)d\nu_n(y)\\\label{eq1} &
+\iint\log(|x-y|^{-1})(1-h(|x-y|))d\mu_n(x)d\nu_n(y).
\end{align}
On one hand, the cartesian product measure $\mu_n\times\nu_n$
tends weakly to $\mu\times\nu$, see \cite[Theorem
2.8]{Bill} or
\cite[Theorem 9.5.9]{Dudley}, and the integrand in the first
integral is lower semi-continuous and lower bounded on $\C$. Hence,
Lemma
\ref{weak-lower-sc} applies (more precisely a version of it on
$\C^2$ which holds true as well). On the other hand, the
second integral has a modulus less than $\epsilon$ uniformly in
$n$.
Consequently,
\begin{align*}
\liminf_{n\to\infty}I(\mu_n,\nu_n) & \geq \iint
\log(|x-y|^{-1})h(|x-y|)d\mu(x)d\nu(y)-\epsilon\\ & =
I(\mu,\nu)-\iint
\log(|x-y|)(1-h(|x-y|))d\mu(x)d\nu(y)-\epsilon.
\end{align*}
The integrand in the last integral is continuous and lower bounded on $\C$. Hence, again by
Lemma \ref{weak-lower-sc}, this integral is less than
$$\liminf_{n\to\infty}\iint
\log(|x-y|)(1-h(|x-y|))d\mu_n(x)d\nu_n(y)\leq\epsilon,$$
which implies
$$\liminf_{n\to\infty}I(\mu_n,\nu_n) \geq I(\mu,\nu)-2\epsilon.$$
Since $\epsilon>0$ is arbitrary, (\ref{I-lsc})
follows.
    \end{proof}

\obeylines \texttt{
Bernhard Beckermann, bbecker@math.univ-lille1.fr
Ana C. Matos, Ana.Matos@univ-lille1.fr
Laboratoire de Math\'ematiques P. Painlev\'e UMR CNRS 8524 - Bat.M2
Universit\'e des Sciences et Technologies Lille
F-59655 Villeneuve d'Ascq Cedex, FRANCE
\medskip
Valeriy Kalyagin, kalia@hse.nnov.ru
Higher School of Economics
Nizhny Novgorod, RUSSIA
\medskip
Franck Wielonsky, wielonsky@cmi.univ-mrs.fr
Laboratoire LATP - UMR CNRS 6632
Universit\'e d'Aix-Marseille
CMI 39 Rue Joliot Curie
F-13453 Marseille Cedex 20, FRANCE
}
\end{document}